\documentstyle[twoside,12pt]{article}
\newtheorem{prop}{Proposition}[subsection]
\newtheorem{opr}[prop]{Definition}
\newtheorem{theo}[prop]{Theorem}

\newtheorem{conj}[prop]{Conjecture}
\newtheorem{rem}[prop]{Remark}
\newtheorem{coro}[prop]{Corollary}
\newtheorem{exam}[prop]{Example}

\newtheorem{lem}[prop]{Lemma}

\title{Generalized Hypergeometric Functions and Rational Curves on Calabi-Yau
Complete Intersections in Toric Varieties}

\author{Victor V. Batyrev  \\
FB Mathematik, Universit\"at-GH-Essen, \\
Universit\"atstra\ss e 3, 45141 Essen, FRG \\
e-mail: matf00@vm.hrz.uni-essen.de \\
and \\
Duco van Straten \\
FB Mathematik, Universit\"at Kaiserslautern \\
Erwin-Schr\"odingerstra\ss e 48, 67663 Kaiserslautern, FRG \\
e-mail: straten@mathematik.uni-kl.de}

\begin{document}

\date{}

\maketitle

\begin{abstract}
We formulate general conjectures about the relationship  between
the A-model connection  on the cohomology of
a $d$-dimensional Calabi-Yau complete intersection
$V$ of $r$ hypersurfaces $V_1, \ldots , V_r$
in a toric variety ${\bf P}_{\Sigma}$
and the system of differential operators annihilating
the special generalized hypergeometric function  $\Phi_0$ depending on
the fan $\Sigma$.
In this context, the
mirror symmetry phenomenon can be interpreted as the twofold
characterization  of the series $\Phi_0$. First,
$\Phi_0$ is  defined by  intersection numbers  of rational curves in
${\bf P}_{\Sigma}$ with the hypersurfaces $V_i$ and their toric
degenerations. Second, $\Phi_0$
is  the power expansion near a boundary point of the moduli space of the
monodromy invariant period of the homolomorphic differential $d$-form on
an another Calabi-Yau $d$-fold $V'$ called the mirror of $V$.
Using this generalized hypergeometric series,
we propose a general construction for mirrors $V'$ of $V$ and canonical
$q$-coordinates on the moduli spaces of Calabi-Yau manifolds.
\end{abstract}

\section{Introduction}

\noindent

In this paper we consider complex projective smooth algebraic varieties
$V$ of dimension $d$ whose canonical bundles ${\cal K}_V$ are trivial,
i.e. ${\cal K}_V \cong {\cal O}_V$, and the Hodge numbers $h^{p,0}(V)$ are
zero unless $p=0$, or $p=d$.  These varieties are called
{\em $d$-dimensional Calabi-Yau varieties}, or {\em Calabi-Yau $d$-folds}.
For each dimension $d \geq 3$, there are many of examples of
topologically different Calabi-Yau $d$-folds which can be constructed
from hypersurfaces and complete intersections in weighted projective
spaces \cite{cand0,cand01,cand02,kim}.

Physicists have discovered a fascinating phenomen for Calabi-Yau manifolds,
so called {\em mirror symmetry} \cite{dixon,greene,lerche,lynker}.
Using the mirror symmetry,
Candelas et al. in  \cite{cand2} have computed  the
coefficients of the $q$-expansion of the Yukawa coupling
for Calabi-Yau hypersurfaces of
degree $5$ in ${\bf P}^4$.
The method of Candelas et al. was applied  to
Calabi-Yau $3$-folds in weighted projective spaces
\cite{font,morrison.picard,klemm1} and
complete intersections in weighted and ordinary projective spaces
\cite{klemm2,lib.teit}.
The  $q$-expansions for Yukawa couplings have been
calculated also for Calabi-Yau hypersurfaces of dimension $d > 3$ in
projective spaces \cite{greene1}.

The interest  of algebraic geometers to
Yukawa couplings is explained by the conjectural
relationship  between the coefficients
of the $q$-expansion of the Yukawa couplings
and the intersection theory on the moduli spaces of
rational curves on Calabi-Yau $d$-folds \cite{greene1,katz1}.
For small values of degrees,
this relationship was verified in \cite{katz2}. However,
the main problem which remains unsolved is to find
a general rigorous mathematical explanation of this relationship
for rational curves of arbitrary degrees on Calabi-Yau $d$-folds.

The purpose of this paper is to show that the mirror symmetry
and the calculation of the  Yukawa couplings for   $d$-dimensional
Calabi-Yau complete intersections in toric
varieties bases essentially on the theory of special generalized
hypergeometric functions. We remark that these hypergeometric
functions satisfy the hypergeometric differential system considered
by Gelfand, Kapranov and Zelevinsky in \cite{gel1}. We propose also a
general method for computing the normalized canonical $q$-coordinates.

The paper is organized as follows:

In Section 2, we give a review of the
calculation of Candelas et al. in \cite{cand2}
of the coefficients $\Gamma_d$ of the $q$-expansion of the normalized
Yukawa coupling
\[ K^{(3)}_q = 5 + \sum_{d \geq 1} \Gamma_d \frac{q^d}{1 -q^d}. \]
The coefficients $\Gamma_d = n_d d^3$  conjecturaly coincide with the
Gromov-Witten invariants  (introduced by D. Morrison
in  \cite{morrison.hodge}) for rational
curves on  quintic hypersurfaces in ${\bf P}^4$.
Our review is greatly influenced by the work of D. Morrison
\cite{morrison.mirror,morrison.picard}, but we want to emphasize on the fact
that the computation of the prediction for the number
of rational curves on  quintic $3$-folds
bases essentially on the properties of the special generalized
hypergeometric series
\[ \Phi_0(z) = \sum_{ n \geq 0} \frac{(5n)!}{(n!)^5}z^n \]
which admits a combinatorial definition in terms of curves on ${\bf P}^4$.

In Section 3, we explain a Hodge-theoretic framework for mirror symmetry
and the ideas due to P. Deligne \cite{deligne1} and
D. Morrison \cite{morrison.comp,morrison.hodge}
The key-point  of this framework is
the existence of a new-type  nilpotent  connection on cohomology of
Calabi-Yau manifolds. Following a suggestion of D. Morrison,  we call it
{\em A-model connection} (see also \cite{witten2}).
The mirror symmetry identifies the
A-model connection on the cohomology of a Calabi-Yau $d$-fold $V$
with the classical Gau\ss -Manin connection on cohomology of
its mirror manifold $V'$.

Section 4 contains a review of the standard  computational technique
based on the recurrent relations satisfied by
coefficients of formal solutions of Picard-Fuchs equations. We use
this techique later
in explicit calculations of $q$-expansions for Yukawa couplings
for some examples of  Calabi-Yau complete intersections in toric varieties.

Section 5 is devoted to complete intersections in ordinary projective
spaces. Using explicit description of
the series $\Phi_0(z)$ for Calabi-Yau complete intersections in projective
spaces, we calculate the $d$-point Yukawa coupling and
propose the explicit construction for mirrors of such Calabi-Yau
$d$-folds for arbitrary dimension $d$.

    In Section 6, we give a general definition of  special
generalized hypergeometric functions and establish the relationships
between these functions and combinatorial properties of
rational curves on toric varieties containing Calabi-Yau complete
intersections. It is easy to see that these generalized
hypergeometric functions form a special subclass of
the generalized hypergeometric functions with {\em resonance} parameters
considered  by Gelfand et al. in \cite{gel1}.
We formulate  general conjectures about the
differential systems and canonial $q$-coordinates
defined by  the generalized hypergeometric series
corresponding to  Calabi-Yau complete
intersections in toric varieties. Using  a
combinatorial interpretation of Calabi-Yau complete
intersections in toric varieties
due to Yu. I. Manin,  we propose an explicit construction  of mirrors.

        In Section 7, we consider in more details
the example  Calabi-Yau hypersurfaces
$V$ of degree $(3,3)$ in ${\bf P}^2 \times {\bf P}^2$.
We use this example for the
illustration of the general computational method we used in Section 8,
where we calculate the $q$-expansions of Yukawa couplings for
some Calabi-Yau complete intersections  in products of projective spaces.
The computations in this section were done by the second author using
an universal computer program based on Maple.

{\bf Acknowledgements. } This work benefited greatly from conversations
with F. Beukers, P. Candelas, B. Greene, D. Morrison, Yu.I. Manin, R. Plesser,
R. Schimmrigk, Yu. Tschinkel and A. Todorov.

We are very grateful to
D. Morrison whose  numerous remarks concerning a preliminary version
of the paper helped us to give precise references on his work, especially  on
the forthcomming papers \cite{greene1,morrison.hodge}

We would like to express
our thaks for
hospitality to the Mathematical Sciences Research Institute where this
work was supported in part by the National Science Foundation
(DMS-9022140), and the DFG (Forschungsschwerpunkt Komplexe
Mannigfaltigkeiten).

\section{Quintics in ${\bf P}^4$}

\noindent

In this section we give a review of the (conjectural) computation of the
Gromov-Witten  invariants $\Gamma_d$ and predictions $n_d$ for
numbers of rational curves of degree $d$ on quintics $V$ in ${\bf P}^5$
due to P. Candelas, X. de la Ossa, P.S. Green,
and L. Parkes \cite{cand2}. The main ingedients of this computations
were considered in papers of D. Morrison
\cite{morrison.mirror,morrison.picard}.
The purpose of this review is
to stress that this computation  needs knowing only
properties of the special generalized hypergeometric
function $\Phi_0(z)$. We begin with the computational algorithm
for computing the coefficients in the $q$-expansion of the
Yukawa coupling and the predictions for number of rational curves.

\newpage

\subsection{The coefficients in the $q$-expansion of the Yukawa coupling}

\noindent

Consider  the series
\[ \Phi_0(z)  = \sum_{n \geq 0} \frac{(5n)!}{(n!)^5}z^n. \]

{\bf Step 1.} If we put $a_n = {\displaystyle \frac{(5n)!}{(n!)^5}}$, then the
numbers $a_n$ satisfy the recurrent relation
\[ (n+1)^4 a_{n+1} = 5(5n+1)(5n+2)(5n+3)(5n+4) a_n. \]
This immediatelly implies that the series $\Phi_0(z)$ is the solution to
the differential equation
\[ {\cal D} \Phi(z) =0 ,\]
where
 \[ {\cal D} = \Theta^4 - 5z(5 \Theta +1)
 (5 \Theta +2)(5 \Theta +3)(5 \Theta +4),\;
 \Theta =  z \frac{\partial}{\partial z}.   \]

One can rewrite the differential operator ${\cal D}$ in powers of
$\Theta$ as follows:
\[ {\cal D} = A_4(z)\Theta^4 + A_3(z)\Theta^3 + \cdots + A_0(z). \]
We denote by $C_i(z)$ the rational function  $A_i(z)/A_4(z)$
$( i =0,\ldots, 3)$.
\bigskip

{\bf Step 2.} Following \cite{morrison.mirror},
define the normalized Yukawa $3$-differential
 as
\[ {\cal W}_{3} = K_z^{(3)} (\frac{d z}{z})^{\otimes 3}, \]
where $K_z^{(3)} = W_3(z)/ \Phi_0^2(z)$ is the $3$-point coupling function
such that $W_3(z)$ satisfies the differential equation
\[ \Theta  W_3(z) = - \frac{1}{2}C_3(z) W_3(z) ({\rm \ref{coup-eq}}) \]
and the normalizing condition $W_3(0) = 5$.

One easily obtains
\[ {\cal W}_{3} = \frac{5}{(1 - 5^5z)\Phi_0^2(z)}
(\frac{d z}{z})^{\otimes 3}. \]
\medskip

{\bf Step 3.} The  equation ${\cal D}\Phi =0$ is a
Picard-Fuchs differential
equation with maximal unipotent monodromy
(in the sense of Morrison \cite{morrison.mirror})
  at $z =0$. Therefore, there exists
a unique solution $\Phi_1(z)$ to ${\cal D}\Phi =0$ such that
$\Phi_1(z) = (\log z)\Phi_0(z) + \Psi(z)$, where $\Psi(z)$ is regular at
$z =0$ and $\Psi(0) =0$. We define the new local coordinate
$q = q(z)$ near the point $z =0$ as
\[ q(z) = \exp \left( \frac{\Phi_1(z)}{\Phi_0(z)} \right)
= z \exp \left( \frac{\Psi(z)}{\Phi_0(z)} \right). \]
The, we rewrite the normalized Yukawa $3$-differential
${\cal W}_{3}$ in the coordinate $q$ as
\[ {\cal W}_{3} = K_q^{(3)} (\frac{d q}{q})^{\otimes 3}. \]
The function $K_q^{(3)}$ is called the Yukawa $3$-point coupling.
This function has the power expansion
\[ K_q^{(3)} = 5 + \sum_{d \geq 1}^{\infty} \frac{n_d d^3 q^d}{1 - q^d}, \]
where $\Gamma_d = n_d d^3$ are conjectured to be the Gromov-Witten invariants
of rational curves of degree $d$ on a
quintic $3$-fold in ${\bf P}^4$ \cite{katz2,morrison.hodge}.   The numbers
$n_d$ are predictions for numbers of rational curves of degree $d$
on quintic $3$-folds.
\bigskip

It is important to remark that in the above algorithm for calculation of
 the numbers $n_d$ one needs to know
only properties of the series $\Phi_0(z)$ and
the normalization condition $W_3(0) = {\rm deg}\, V = 5$ for $W_3(z)$,
i.e., one does not
need to know anything about mirrors of quintics.

\subsection{Philosophy of mirrors and the series $\Phi_0(z)$}

\noindent

The central role in the computation of Candelas et al. in \cite{cand2}
is played by  the orbifold construction of mirrors for quintics in ${\bf P}^4$
\cite{greene}.
In \cite{batyrev.dual}, this construction of mirrors was generalized for
hypersurfaces in toric Fano varieties with Gorenstein singularities.

In the above algorithm, we have shown that one can forget about
mirrors. However, the philosophy of mirrors proves to be
very helpful. For quintic $3$-folds this philosophy appears as the following
twofold  interpretation  of the series $\Phi_0(z)$.
\medskip

{\bf The first interpretation:} We compute the coefficients $a_n$ of the
power series $\Phi_0(z)$ using combinatorial
properties of curves $C \subset {\bf P}^4$ of degree $n$.

Notice that
any such a curve
$C$ meets a generic quintic $V$ at $5n$ distinct points
$ p_1,$ $\ldots,$ $p_{5n}$. There exists a degeneration of $V$ into a union of
$5$ hyperplanes $H_1 \cup \cdots \cup H_5$. Every such a hyperplane $H_i$
intersects $C$ at $n$ points $p_{i_1},$ $\ldots,$ $p_{i_n}$ which can be
considered as deformations of a subset of $n$ points from the set $\{ p_1,$
 $\ldots,$ $ p_{5n} \}$.
It remains to remark that there exists exactly $(5n)!/(n!)^5$ ways to divide
$\{ p_1,$ $\ldots,$ $p_{5n} \}$ into $5$ copies of
$n$-element disjoint subsets.
\medskip

{\bf The second interpretation: } We find the coefficients $a_n$ from
an integral representation of $\Phi_0(z)$.

Let ${\bf T} \cong ({\bf C}^*)^4$ be
the $4$-dimensional algebraic torus with coordinate functions
$X_1$, $X_2$, $X_3$, $X_4$. Take  the Laurent polynomial
\[f(u,X) = 1 - (u_1X_1 + u_2X_2 + u_3X_3 + u_4 X_4 +
u_5(X_1 X_2 X_3 X_4)^{-1})  \]
in variables $X_1, X_2, X_3, X_4$, where the coefficients
$u_1, \ldots, u_5$ are considered as independent parameters. Let
$z = u_1 u_2 u_3 u_4 u_5$.

\begin{prop}
\[ \Phi_0(u_1 \cdots u_5) = \Phi_0(z) = \frac{1}{(2 \pi \sqrt{-1})^4}
\int_{\mid X_i \mid =1} \frac{1}{f(u,X)}
\frac{dX_1}{X_1} \wedge \frac{dX_2}{X_2} \wedge
\frac{dX_3}{X_3} \wedge \frac{dX_4}{X_4}. \]
\end{prop}

{\bf Proof. } One has
\[ \frac{1}{f(u,X)} = \sum_{n \geq 0} (u_1X_1 + u_2X_2 + u_3X_3 + u_4 X_4 +
u_5(X_1 X_2 X_3 X_4)^{-1})^n \]
 \[ = \sum_{m \in {\bf Z}^4} c_m(u) X^m. \]
It is straightforward to see that $c_0(u) = \Phi_0(u_1 \cdots u_5)$.
Now the statement follows from the Cauchy residue formula. \hfill $\Box$
\medskip

The second interpretation of $\Phi_0(z)$ implicitly uses mirrors of
quintics, since zeros of $f(u,X)$ define the affine
Calabi-Yau $3$-fold $Z_f$ in ${\bf T}$ whose smooth Calabi-Yau
compactification  is mirror symmetric with respect to
quintic $3$-folds (see \cite{batyrev.dual}).
Moreover, the holomorphic $3$-form $\omega$
on $Z_f$
such that $\omega(z)$ extends regularly to a smooth compactification
of $Z_f$ depends only on $z$, i.e., only on the product $u_1 \cdots u_5$.
This $3$-form   can be written  as
\[ \omega(z) = \frac{1}{(2\pi \sqrt{-1})^4}{\rm Res}
\frac{1}{f(u,X)} \frac{dX_1}{X_1} \wedge \frac{dX_2}{X_2} \wedge
\frac{dX_3}{X_3} \wedge \frac{dX_4}{X_4}. \]
This shows that $\Phi_0(z)$ is exactly the monodromy
invariant period of the $3$-form
$\omega(z)$ near $z =0$.

\begin{prop}
The differential $3$-form $\omega(z)$ satisfies the same Picard-Fuchs
differential equation ${\cal D} \Phi = 0$ as for the series $\Phi_0(z)$.
In particular, all periods of $\omega(z)$ satisfy the Picard-Fuch
differential equation with the operator
\[ \Theta^4 - 5z(5 \Theta +1)
 (5 \Theta +2)(5 \Theta +3)(5 \Theta +4). \]
\end{prop}

{\bf Proof.} In order to prove  the statement, it is sufficient  to check that
\[ \left({\cal D} \frac{1}{f(u,X)} \right) \frac{dX_1}{X_1}
\wedge \frac{dX_2}{X_2} \wedge
\frac{dX_3}{X_3} \wedge \frac{dX_4}{X_4}  \]
is a differential of a rational $3$-form on ${\bf T} \setminus Z_f$ .
The latter follows from standard arguments using reduction
by the Jacobian ideal $J_f$ (see
\cite{batyrev.var}). \hfill $\Box$

\subsection{A-model connection}

\noindent

The Yukawa coupling can be described by a nilpotent connection $\nabla_A$
on the cohomology of quintic $3$-fold $V$
\[ \nabla_A \; : \; H^*(V, {\bf C}) \rightarrow H^*(V, {\bf C}) \otimes
{\bf C} \langle  \frac{dz}{z} \rangle. \]
This connection is homogeneous of degree $2$, i.e.,
\[ \nabla_A\; : \; H^i (V, {\bf C}) \rightarrow H^{ i +2}(V, {\bf C})\otimes
{\bf C} \langle  \frac{dz}{z} \rangle, \]
and hence $\nabla_A$ vanishes on $H^3(V, {\bf C})$. By this reason, we consider
only  the cohomology subring
\[ H^{2*}(V, {\bf Z}) = \oplus_{i =0}^3  H^{2i}(V, {\bf Z}) \subset
H^*(V, {\bf Z})\]
of even-dimensional classes on a quintic $3$-fold  $V$
(${\rm rk}\,H^{2i} (V, {\bf Z}) =1$).
Let $\eta_i$ be the positive generator of  $H^{2i} (V, {\bf Z})$. Then,
in the basis $\eta_0, \eta_1, \eta_2, \eta_3$,
 the multiplication by $\eta_1$ is the endomorphism of
$H^{2*}(V, {\bf Z})$ having as
matrix
\[  \Lambda = \left( \begin{array}{cccc} 0 & 1 & 0 & 0 \\
0 & 0 & 5 & 0 \\ 0 & 0 & 0 &
1 \\
0 & 0  & 0  &  0 \end{array} \right). \]

Following \cite{cere} and \cite{morrison.hodge},
we define the $1$-parameter connection on
$H^{2*}(V, {\bf C}) \otimes {\bf C} [[q]]$ considered as a trivial
bundle over ${\rm Spec}\, {\bf C}[[q]]$ as follows
\[ \left( \begin{array}{c}  \nabla_A \eta_0 \\
\nabla_A \eta_1 \\ \nabla_A \eta_2 \\ \nabla_A \eta_3 \end{array} \right)
= \left( \begin{array}{cccc} 0 & \frac{dq}{q} & 0 & 0 \\
0 & 0 &  K_q^{(3)}\frac{dq}{q} & 0 \\ 0 & 0 & 0 &
\frac{dq}{q} \\
0 & 0 & 0 & 0 \end{array} \right)
\left( \begin{array}{c}   \eta_0 \\ \eta_1 \\  \eta_2 \\
\eta_3 \end{array} \right). \]

\begin{equation}
 K(q) =   \left( \begin{array}{cccc} 0 & \frac{dq}{q} & 0 & 0 \\
0 & 0 &  K_q^{(3)}\frac{dq}{q} & 0 \\ 0 & 0 & 0 &
\frac{dq}{q} \\
0 & 0 & 0 & 0 \end{array} \right)  \label{matrix}
\end{equation}
can be considered as the deformation of the matrix $\Lambda$ such that
\[ \Lambda = {\rm Res}\mid_{q =0}\, K(q). \]
\bigskip

The mirror philosophy shows that the matrix (\ref{matrix})
can be identified with the matrix of the classical Gau\ss -Manin
connection on the $4$-dimensional cohomology space
$H^3(\hat{Z}_f, {\bf C})$ in a special
symplectic basis. We notice that the quotients ${ F}^i/{ F}^{i+1}$
of the Hodge filtration
\[ H^3(\hat{Z}_f, {\bf C}) = { F}^0 \supset { F}^1 \supset
{F}^2 \supset { F}^3 \supset {F}^4 = 0 \]
are $1$-dimensional. There is also the monodromy filtration on the
homology $H_3(\hat{Z}_f, {\bf Z})$
\[ 0 = { W}_{-4} \subset { W}_{-3} \subset { W}_{-2} \subset
{ W}_{-1} \subset { W}_{0} = H_3(\hat{Z}_f, {\bf Z}) \]
such that ${ W}_i/{ W}_{i-1}$ are also $1$-dimensional.
We  choose the symplectic basis $\gamma_0,\gamma_1, \gamma_2,\gamma_3$ in
$H_3(\hat{Z}_f, {\bf Z})$ in such a way that
$\{ \gamma_0, \ldots, \gamma_i \}$ form a ${\bf Z}$-basis of ${W}_i$.
We choose also the basis $\omega_0, \omega_1, \omega_2, \omega_3$ of
$H^3(\hat{Z}_f, {\bf C})$ such that
$\{ \omega_0, \ldots, \omega_i \}$ form a ${\bf C}$-basis of ${F}^{3-i}$
and
\[ p_{ij} = \int_{\gamma_j} \omega_i = \delta_{ij}\;\;
{\rm for}\;i \geq j. \]
So the period matrix $\Pi = (p_{ij})$ has the form
\cite{greene1,morrison.hodge}
\[ \Pi = \left( \begin{array}{cccc} 1 & p_{12} & p_{13} & p_{14} \\
0 & 1 & p_{23} & p_{34} \\ 0 & 0 & 1 & p_{34} \\
0 & 0 & 0 & 1 \end{array} \right). \]
Notice that all coefficients $p_{ij}$ $(i < j)$ are multivalued functions
of $z$ near $z = 0$.

Applying the Griffiths transversality property, we obtain
that  the
Gau\ss--Manin connection in the $z$-coordinate has the form
\[ \left( \begin{array}{c}  \nabla \omega_0 \\
\nabla \omega_1 \\ \nabla \omega_2 \\ \nabla \omega_3 \end{array} \right)
= \left( \begin{array}{cccc} 0 & (\Theta p_{12})\frac{dz}{z} & 0 & 0 \\
0 & 0 & (\Theta p_{23})\frac{dz}{z} & 0 \\ 0 & 0 & 0 &
(\Theta p_{34})\frac{dz}{z} \\
0 & 0 & 0 & 0 \end{array} \right)
\left( \begin{array}{c}   \omega_0 \\ \omega_1 \\  \omega_2 \\
\omega_3 \end{array} \right), \]
where $\Theta p_{i,i+1}$ are single valued functions.

Then the Yukawa $3$-differential is simply the tensor product
\[ {\cal W}_3 = K_z^{(3)} (\frac{dz}{z})^{\otimes 3} =
(\Theta p_{12})\frac{dz}{z} \otimes (\Theta p_{23})\frac{dz}{z} \otimes
(\Theta p_{34})\frac{dz}{z}. \]
\medskip

By Griffiths transversality, one has $\omega_0 \wedge \omega_2 = 0$,
i.e. we can assume that $p_{23} = p_{12}$. The
differential form $w_0$ can be defined as $\omega / \Phi_0(z)$.
Moreover, $p_{12} = \Phi_1(z) / \Phi_0(z)$. In the new
coordinate $q$,  we have  $p_{12} = \log q$.
Then the Gau\ss -Manin connection can be rewritten as
\[ \left( \begin{array}{c}  \nabla \omega_0 \\
\nabla \omega_1 \\ \nabla \omega_2 \\ \nabla \omega_3 \end{array} \right)
= \left( \begin{array}{cccc} 0 & \frac{dq}{q} & 0 & 0 \\
0 & 0 &  K_q^{(3)}\frac{dq}{q} & 0 \\ 0 & 0 & 0 &
\frac{dq}{q} \\
0 & 0 & 0 & 0 \end{array} \right)
\left( \begin{array}{c}   \omega_0 \\ \omega_1 \\  \omega_2 \\
\omega_3 \end{array} \right). \]
The role played by the  $1$-form
$dq/q$ in the connection matrix was first observed by
 P. Deligne \cite{deligne1}.

\subsection{The $q$-coordinate and the Yukawa coupling}

Since the coordinate $q$ was defined intrinsically as
the ratio $\Phi_1(z)/\Phi_0(z)$
of two solutions of the differential equation ${\cal D} \Phi =0$, it is
natural to ask about the form of the differential operator ${\cal D}$ in
the new coordinate $q$. Denote by $\Xi$ the
differential operator ${\displaystyle q \frac{\partial}{\partial q}}$.

\begin{prop}
The differential $3$-form $\omega_0$ satisfies the Picard-Fuchs
differential  equation with the differential operator
\[ \Xi^4 + c_3(q) \Xi^3 + c_2(q)  \Xi^2,  \]
where
\[ c_3(q) = - 2\frac{\Xi K_q^{(3)}}{K_q^{(3)}}, \;
c_2(q) =  \frac{\Xi K_q^{(3)}}{(K_q^{(3)})^2} -
\frac{\Xi^2 K_q^{(3)}}{K_q^{(3)}}. \]
\end{prop}

{\bf Proof. } By properties of the nilpotent connection, one has
\[ \Xi \omega_0 = \omega_1, \; \Xi \omega_1 = K_q^{(3)} \omega_2, \;
\Xi \omega_2 = \omega_3, \; \Xi \omega_3 = 0. \]
So
\[ \Xi^4 \omega_0 = \Xi^2 K_q^{(3)} \omega_2 = \Xi( (\Xi K_q^{(3)}) \omega_2 +
K_q^{(3)} \omega_3  )  \]
\[ = (\Xi^2 K_q^{(3)}) \omega_2 +  2 (\Xi K_q^{(3)}) \omega_3. \]
On the other hand,
\[ \omega_2 = \frac{1}{K_q^{(3)}} \Xi^2 \omega_0, \]
\[ \omega_3 = \Xi(\frac{1}{K_q^{(3)}} \Xi^2 \omega_0)  =
- \frac{\Xi K_q^{(3)}}{(K_q^{(3)})^2} \Xi^2 \omega + \frac{1}{K_q^{(3)}}\Xi^3
\omega_0. \]
\hfill $\Box$

\begin{rem}
{\rm  The differential equation for $\omega_0$ can be written also as
\[ \Xi^2 (K_q^{(3)})^{-1} \Xi^2 \omega_0 = 0. \]
In this form this equation first arose in \cite{ferrara.louis}.  }
\end{rem}

The differential operator ${\cal D}$ which annihilates the
function $\Phi_0(z)$ defines the connection in the basis
 $\omega$, $\Theta \omega$, $\Theta^2 \omega$,
$\Theta^3 \omega$ of $H^3(\hat{Z}_f, {\bf C})$ :
\[ \left( \begin{array}{c}  \nabla\, \omega \\
\nabla\, \Theta \omega \\ \nabla\, \Theta^2 \omega
\\ \nabla\, \Theta^3 \omega
\end{array} \right)
= \left( \begin{array}{cccc} 0 & \frac{dz}{z} & 0 & 0 \\
0 & 0 & \frac{dz}{z} & 0 \\ 0 & 0 & 0 &
\frac{dz}{z} \\
-C_0(z)\frac{dz}{z} & -C_1(z)\frac{dz}{z}  & -C_2(z)\frac{dz}{z} &
-C_3(z)\frac{dz}{z} \end{array} \right)
\left( \begin{array}{c}   \omega \\
 \Theta \omega \\  \Theta^2 \omega \\ \Theta^3 \omega
\end{array} \right). \]
The basis $\omega$, $\Theta \omega$, $\Theta^2 \omega$,
$\Theta^3 \omega$ is also
compatible with the Hodge filtration in $H^3(\hat{Z}_f, {\bf C})$.
Thus there exist a matrix
\[ R = \left( \begin{array}{cccc} r_{11} & r_{12} & r_{13} & r_{14} \\
0 & r_{22} & r_{23} & r_{34} \\ 0 & 0 & r_{33} & r_{34} \\
0 & 0 & 0 & r_{44} \end{array} \right) \]
such that
\[ \left( \begin{array}{c}   \omega \\
 \Theta \omega \\  \Theta^2 \omega \\ \Theta^3 \omega
\end{array} \right) = R
\left( \begin{array}{c}   \omega_0 \\ \omega_1 \\  \omega_2 \\
\omega_3 \end{array} \right). \]
It is easy to see that
\[ r_{11} = \Phi_0(z), \; r_{22} = \Phi_0(z) (\Theta p_{12}), \;
r_{33} = \Phi_0(z) (\Theta p_{12})(\Theta p_{23}), \; \]
\[ r_{44} = \Phi_0(z) (\Theta p_{12})(\Theta p_{23})(\Theta p_{34}). \]
\bigskip

\section{Quantum variations of Hodge structure on Calabi-Yau
manifolds}

\subsection{A-model connection  and rational curves}

\noindent

A general appoach to the definition of a new
connection on cohomology
of algebraic and symplectic manifolds $V$
was proposed by Witten \cite{witten}.
The construction of Witten bases on the interpretation of
third partial derivatives
\[ \frac{\partial^3 }{\partial z_i \partial z_j \partial z_k} P(z) \]
of a function $P(z)$ on the cohomology space $H^*(V, {\bf C})$ as structure
constants of a commutative associative algebra. The
function $P(z)$ is defined via the intersection theory on the
moduli spaces of mappings of Riemann surfaces $S$ to $V$. Using Poincare
duality, one obtains the structure coefficients of the connection on
$H^*(V, {\bf C})$.

We consider a specialization of the general construction to the case
when $V$ is a Calabi-Yau $3$-fold. We put $n = {\rm dim}\, H^2 (V, {\bf C}) =
{\rm dim}\, H^4(V, {\bf C})$. Let $\eta_0$ be a generator of
$H^0(V, {\bf Z})$, $\eta_1, \ldots, \eta_n$ a ${\bf Z}$-basis of
$H^2(V, {\bf Z})$, $\zeta_1, \ldots, \zeta_n$ the dual ${\bf Z}$-basis of
$H^4(V, {\bf Z})$ $( \langle \eta_i, \zeta_j \rangle = \delta_{ij})$, and
$\zeta_0$ the dual to $\eta_0$ generator of
$H^6(V, {\bf Z})$. We can always assume that the
cohomology classes   $\eta_1, \ldots, \eta_n$ are contained in the
closed K\"ahler cone of $V$.

\begin{opr}
{\rm Let $R = {\bf C} [[ q_1, \ldots, q_n ]]$ be the ring of formal
power series
in $n$ independent variables. We denote by $H(V)$ the scalar extension
\[ (\oplus_{i =0}^3 H^{2i}(V, {\bf C})) \otimes_{\bf C} R. \]
We consider a flat nilpotent holomorphic connection
\[ \nabla_A\; : H(V) \rightarrow H(V) \otimes \Omega^1_R(\log q) \]
defined by  the following formulas \cite{cere,morrison.hodge}

\[  \nabla_A \eta_0 = \sum_{i =1}^n \eta_i \otimes \frac{dq_i}{q_i}; \]

\[ \nabla_A \eta_k = \sum_{i =1}^n \sum_{j =1}^n K_{ijk} \; \zeta_j \otimes
\frac{d q_i}{q_i},\; k =1, \ldots, n;  \]

\[ \nabla_A \zeta_j = \zeta_0 \; \frac{d q_j}{q_j} , \; j =1, \ldots, n; \]
\[ \nabla_A \zeta_0 = 0. \]
The coefficients $K_{ijk}$ are power series in $q_1, \ldots, q_n$
defined by rational curves $C$ on $V$, i.e., morphisms
$f \; :\; {\bf P}^1 \rightarrow V$ as follows
\[ K_{ijk} = \langle \eta_i, \eta_j, \eta_k \rangle +
\sum_{\scriptstyle \begin{array}{c} {\scriptstyle C \subset V} \\
{\scriptstyle \lbrack C \rbrack \neq 0} \end{array} } n_{\lbrack C \rbrack}
\langle C, \eta_i \rangle \langle C, \eta_j \rangle
\langle C, \eta_k \rangle \frac{q^{\lbrack C \rbrack} }
{ 1 - q^{\lbrack C \rbrack}}, \]
where $q^{\lbrack C \rbrack} = q_1^{c_1} \cdots  q_n^{c_n}$ $( c_i =
\langle C, \eta_i \rangle)$. The integer
\[ \Gamma_{\lbrack C \rbrack}(\eta_i,\eta_j,\eta_k) =
n_{\lbrack C \rbrack} \langle C, \eta_i \rangle \langle C, \eta_j \rangle
\langle C, \eta_k \rangle \]
is  called {\em the Gromov-Witten invariant}
\cite{katz2,morrison.hodge} {\em of
the class $\lbrack C \rbrack$.} If the classes $\eta_i$, $\eta_j$, and
$\eta_k$ are represented by effective divisors $D_i$, $D_j$, and $D_k$ on
$V$, then $\Gamma_{\lbrack C \rbrack}(\eta_i,\eta_j,\eta_k)$ is
the number of pseudo holomorphic immersions $\imath \; : \; {\bf P}^1
\rightarrow V$ such that $\lbrack \imath({\bf P}^1) \rbrack =
\lbrack C \rbrack$ and $\imath(0) \in D_i$, $\imath(1) \in D_j$,
$\imath(\infty) \in D_k$ for sufficiently general almost complex structure
on $V$. In particular, the number $n_{\lbrack C \rbrack}$ which predicts
the number of rational curves $C \subset V$ with the fixed
class $\lbrack C \rbrack$ is
always  non-negative.

The connection $\nabla_A$ will be called the {\em A-model connection}.
The  connection $\nabla_A$ defines on $H(V)$  a
variation of Hodge structure of type $(1,n,n,1)$. We call this
variation {\em the quantum  variation of Hodge structure on $V$. } }
\end{opr}

\begin{rem}
{\rm The Picard-Fuch differential
system satisfied by $\eta_0$ was considered in details in  \cite{cere}.}
\end{rem}

One immediatelly obtains:

\begin{prop}
Let $\eta = l_1 \eta_1 + \cdots + l_n \eta_n \in H^2(V, {\bf Z})$
be a class of an ample divisor on $V$. Define the $1$-parameter
connection with the new coordinate $q$ by putting
$q_1 =q^{l_1}\,  \ldots, q_n = q^{l_n}$.  Then the connection $\nabla_A$ on
$H(V)$ induces the connection
\[ \nabla_q\; : \; (\oplus_{i =0}^3 H^{2i}(V, {\bf C}[[q]]))
 \rightarrow
(\oplus_{i =0}^3 H^{2i}(V, {\bf C}[[q]])) \otimes_{\bf C}
\Omega^1_{{\bf C}[[q]]}(\log q). \]
In particular, the residue of the connection
operator $\nabla_q$ at $q = 0$ is the
Lefschetz operator $L_{\eta} \;:  H^{2i}(V, {\bf C}) \rightarrow
H^{2i +2}(V, {\bf C})$, and
\[ \langle (\nabla_q)^3 \eta_0, \zeta_0 \rangle =
\left( \langle \eta, \eta, \eta \rangle +
\sum_{d > 0} n_d \frac{d^3 q^d}{1 - q^d}  \right) \left( \frac{dq}{q}
\right)^{\otimes 3} \]
where
\[ n_d = \sum_{ \langle C, \eta \rangle = d } n_{\lbrack C \rbrack}. \]
\end{prop}

\begin{coro}
The connection $\nabla_q$ defines  a Picard-Fuchs differential operator of
order $4$ annihilating $\eta_0$.
\end{coro}

\subsection{The Gau\ss -Manin connection for mirrors}

\noindent

Let $W$ be a Calabi-Yau $3$-fold such that ${\rm dim }\,H^3(W, {\bf C}) =
2n +2$. Assume that we are given a variation $W_z$ of complex structure on $W$
near a boundary point $p$ of the $n$-dimensional moduli space ${\cal M}_W$
of complex structures on $W$ in  holomorphic coordinates
$z_1, \ldots, z_n$ near $p$ such that $p = (0,\ldots, 0)$.

\begin{opr}
{\rm The family $W_z$ is said to have the {\em  maximal unipotent
monodromy at $z = 0$} if the weight filtration
\[ 0 = W_{-1} \subset W_0 \subset W_1 \subset W_{2} \subset W_3 = H^3(W_z,
{\bf C}) \]
defined by $N$ is orthogonal to the Hodge filtration $\{ F^i \}$, i.e.,
\[ H^3(W_z, {\bf C}) = W^{\perp}_{i} \oplus F^{3-i} \; i =0, \ldots , 3.\]}
\end{opr}
(This is essentially the same definition given in
\cite{morrison.comp,morrison.hodge}.)

Choose a symplectic basis
\[ \gamma_0, \gamma_1 , \ldots, \gamma_n, \delta_1, \ldots \delta_n,
\delta_0 \]
of  $H_3(W_z, {\bf Z})$ in such a way that
$\gamma_0$ generates $W_0$, $\gamma_0, \gamma_1, \ldots , \gamma_n$ is a
${\bf Z}$-basis of $W_1$,
\[ \gamma_0, \gamma_1 , \ldots, \gamma_n, \delta_1,
\ldots \delta_n \]
is a ${\bf Z}$-basis of $W_2$.
Then we choose a symplectic basis in $H^3(W_z, {\bf C})$ :
\[ \omega_0, \omega_1, \ldots, \omega_n, \nu_1, \ldots, \nu_n, \nu_0 \]
such that $\omega_0$ generates $F^3$,
$\omega_0, \omega_1, \ldots, \omega_n$ is the basis of
$F^2$, $\omega_0, \omega_1, \ldots, \omega_n, \nu_1, \ldots, \nu_n$ is
the basis of $F^1$ such that
\[ \langle \omega_i, \gamma_i \rangle =
\langle \nu_i, \delta_i \rangle = 1,\; i =0,1, \ldots, n \; \]

\[ \langle \omega_i, \gamma_0 \rangle
 = \langle \nu_j, \gamma_0 \rangle =
\langle \nu_i, \gamma_j \rangle  = \langle \nu_0, \gamma_j \rangle =
\langle \nu_0, \delta_j \rangle = 0, \;
 i =1, \ldots, n, \; j =0, \ldots, n. \]
The choice of the basis of $H^3(W_z, {\bf C})$ defines the splitting into
the direct sum
\[ H^3(W_z, {\bf C}) = H^{3,0} \oplus H^{2,1} \oplus H^{1,2} \oplus H^{0,3} \]
such that all direct summand acquire {\em canonical integral structures}.
By Griffiths transversality property, the Gau\ss -Manin connection $\nabla$
sends $H^{3-i,i}$ to $H^{3-i-1, i+1} \otimes \Omega^1(\log z)$.
It is an observation of Deligne \cite{deligne1} that the weight and
Hodge filtrations define a variation of mixed Hodge structure
(VMHS).

Two Calabi-Yau $3$-folds $V$ and $W$ are called mirror symmetric if
the quantum variation of Hodge structure for $V$ is isomorphic to the
classical VMHS
for $W$. In this case the $q$-coordinates near $p$ up to constants are
defined by the formula \cite{morrison.comp}

\[ q_i = \exp (2\pi \sqrt{-1}) \int_{\gamma_i} \omega_0. \]
\bigskip

\section{Picard-Fuchs equations}

\noindent

In this section we recall standard facts about Picard-Fuch
differential equations which we use in computations of
Yukawa $d$-point functions and predictions for numbers of  rational curves
on Calabi-Yau manifolds.

\subsection{Recurrent relations and differential equations}

\noindent

Let $a_n$ $n =0,1,2, \ldots$ be an infinite sequence of complex numbers.
For our purposes, it will be more convenient to define $a_n$ for all
integers  $n \in {\bf Z}$ by putting $a_n =0$ for $n < 0$.
We define the generating function for
the sequence $\{ a_i \}$ as the formal power series
\[ \Phi(z) = \sum_{i \geq 0}^{\infty} a_i z^i \in {\bf C}[[ z
]]. \]

Consider two differential operators acting on ${\bf C} [[ z ]]$:
\[ \Theta \; :\; f \mapsto  z \frac{\partial}{\partial z} f,   \]
\[ z \; : \; f \mapsto z \cdot f \]
satisfying the relation
\begin{equation}
 \lbrack \Theta,z \rbrack = \Theta \circ z - z \circ \Theta = z.
\label{commut}
\end{equation}

These operators generate the algebra
${\bf D} ={\bf C} \lbrack z, \Theta \rbrack$ of "logarithmic"
differential operators which are polynomials in non-commuting
operators $\Theta$ and $z$
\bigskip

Fix a positive integer $d$. Assume that there exist
$m+1$ $(m \geq 1)$ polynomials
\[ P_{0}(y), \ldots, P_m(y) \in {\bf C}\lbrack  y \rbrack \]
of degree $d+1$ such that that for every $n \in {\bf Z}$
the numbers $\{a_i \}$ satisfy the recurrent relation:
\begin{equation}
  P_0(n)a_{n} + P_1(n+1)a_{n +1} + \cdots + P_m(n+m) a_{n+m} = 0.
\label{recurrent}
\end{equation}
(Here we consider $y$ as a new complex variable having no connection to
our previous variable $z$.)   Then $\Phi(z)$ is a formal solution of the
linear differential equation
\[ {\cal D} \Phi(z) = 0\]
with the differential operator
\begin{equation}
 {\cal D} =  z^m P_0(\Theta) + z^{m-1}P_1(\Theta) + \cdots +
P_m(\Theta).
\label{diff1}
\end{equation}
This differential equation of order $d+1$  can be rewritten in
powers of $\Theta$ as
\begin{equation}
{\cal D} = A_{d+1}(z)\Theta^{d+1} + \cdots + A_1(z) \Theta + A_0(z),
\label{diff2}
\end{equation}
where $A_i$ are some polynomials in $z$.
It is easy to check the following:

\begin{prop}
A power series $\Phi(z)$ is a formal solution to a differential
equation ${\cal D}\Phi(z) = 0$  of order $d+1$ for some element ${\cal D} \in
{\bf D}$ if and only if the coefficients $\{ a_i \}$ satisfy a recurrent
relation as in {\rm (\ref{recurrent})} for some polynomials
$ P_{0}(y), \ldots, P_m(y)$ of degree $d+1$.
\end{prop}

\subsection{Picard-Fuchs operators}

\noindent

Recall that a differential operator ${\cal D}$ as in (\ref{diff2}) is called
{\em a Picard-Fuchs operator at point $z =0$} if $A_{d+1}(0) \neq 0$.
Moreover, the condition that solutions of the
Picard-Fuchs equatios ${\cal D}\Phi$ have {\em
maximal unipotent monodromy
at $z =0$} \cite{morrison.picard} is equivalent to
$A_i(0) = 0$ for $ i = 0, \ldots, d$.
These properties of the operator ${\cal D}$  can be
reformulated in terms of the properties
of the polynomial $P_m(y)$ in (\ref{recurrent})
as follows:
\medskip

{\em A differential operator ${\cal D}$ is a Picard-Fuchs operator if and
only if the polynomial $P_m(y)$ has degree $d+1$, i.e., its leading
coefficient is nonzero. Moreover,
solutions of the equations ${\cal D}$ have maximal unipotent monodromy
at $z =0$ if and only if
the polynomial $P_m(y)$ equals $cy^k$ for some nonzero constant $c$. }
\medskip

Picard-Fuchs operators having the  maximal unipotent monodromy at $z =0$
will be
objects of our main interest. Therefore, we introduce the following
definition:

\begin{opr}
{\em A Picard-Fuchs operator ${\cal D}$ with the  maximal unipotent
monodromy will be called {\em a MU-operator}. We will  always assume
that the corresponding polynomial $P_m(y)$ in (\ref{recurrent}) for
any $MU$-operator ${\cal D}$ is $y^k$, i.e., $c =1$.}
\label{MU}
\end{opr}

The fundamental property of $MU$-operators is the following one:

\begin{theo}
If ${\cal D}$ is $MU$-operator, then the subspace in ${\bf C}[[z]]$
of solutions of the linear differential equation
\[ {\cal D}\Phi(z) = 0 \]
has dimension $1$. Moreover, every  solution is defined uniquely by
the value $\Phi(0) = a_0$.
\end{theo}

{\bf Proof. } If we have chosen a value of $a_0$, the all coefficients
$a_i$ for $i \geq 0$ are uniquely defined from the reccurent relation
(\ref{recurrent}). ( We remind, that we put $a_i =0$ for $ i < 0$.) $\Box$

\begin{opr}
{\rm Let ${\cal D}$ be a $MU$-operator. Then the power series
solution $\Phi_0(z)$ of
the equation ${\cal D}\Phi(z) = 0$ normalized by the condition
$\Phi_0(0) =1$ will be called the {\em socle-solution}. }
\end{opr}

\subsection{Logarithmic solutions and the $q$-coordinate}

\noindent

Let ${\cal D}$ be a $MU$-operator of order ${d+1}$. Putting
$C_i(z) = A_i(z)/A_{d+1}(z)$  we can define another
differential operator
\[ {\cal P} = {\cal P}(\Theta) =  \sum_{i =0}^{{d+1}} C_i(z)\Theta^i  \]
which is also a $MU$-operator of order ${d+1}$, where $C_i(z)$ are
rational functions in $z$, and $C_{d+1}(z) \equiv 1$. Assume that we have
a formal regular solution
\[ \Phi(z) = \sum_{i =0}^{\infty} a_n z^n. \]
Consider a formal polynomial extension
\[M_{z} = {\bf C} [[  z ]] \lbrack \log z \rbrack, \]
where $\log z$ is considered as a new transcendent variable. We can define
the structure of a left ${\bf D}$-module on $M_{z}$ putting
by definition $\Theta \log z = 1$. In fact, $M_{z}$ will be a module
over the larger algebra ${\bf D}_z$ containing the new operator $Log\,z$ such
that
\[ z \circ(\Theta \circ Logz\,) =  (\Theta \circ Log\,z) \circ z = 1, \]
\[ \Theta \circ Log\,z -  Log\,z \circ \Theta  =1, \]
and $Log\,z$ acts on $M_z$ by multiplication on $\log z$.

\begin{prop}
Let ${\cal P} = \sum_{i =0}^{{d+1}} C_i(z)\Theta^i$ be any operator in
${\bf D}$. Then
\[ \lbrack {\cal P}, Log\, z \rbrack =
\sum_{i =1}^{{d+1}} iC_i(z)\Theta^{i-1} =
{\cal P}_{\Theta}', \]
where ${\cal P}_{\Theta}'$ is a formal derivative of ${\cal P}$ with respect
to $\Theta$.
\end{prop}

{\bf Proof.} The statement follows from relation
\[ \Theta^i \circ Log\, z - Log\, z \circ \Theta^i = i \Theta^{i-1} \]
which can be proved by induction. \begin{flushright}$\Box$\end{flushright}
\bigskip

Assume that we want to find a element $\Phi_1(z)$ in  $M_z$ such that
${\cal P} \Phi_1(z) = 0$ and  $\Phi_1(z)$ has  form
\[ \Phi_1(z) = \log z \cdot \Phi_0(z) + \Psi(z) \]
where $\Psi(z)$ is an element of ${\bf C} [[  z ]]$, and  $\Psi(0) =0$.

\begin{prop}
The element $\Psi(z)$ satisfies the linear non-homogeneous
differential equation
\begin{equation}
 {\cal P}_{\Theta}' \Phi_0(z) +
{\cal P} \Psi(z) = 0,
\label{exp-diff}
\end{equation}
or, formally,
\[ \Psi(z) = - {\cal P}^{-1} {\cal P}_{\Theta}'\Phi_0(z) = \partial_{\Theta}
\log {\cal P} \cdot \Phi_0(z). \]
\label{solut}
\end{prop}

{\bf Proof.}  Since $\Phi_0$ and $\Phi_1$ are solutions, we obtain
\[ 0 = {\cal P}\Phi_1 = {\cal P} \log z \Phi_0 + {\cal P} \Psi = \]
\[ = ( Log \,z \circ{\cal P} + \lbrack {\cal P}, Log\, z \rbrack)\Phi_0 +
{\cal P}\Psi
 = \lbrack {\cal P}, Log\, z  \rbrack \circ \Phi_0(z) +
{\cal P}\Psi = {\cal P}_{\Theta}' \Phi_0 + {\cal P}\Psi_0 = 0.\;  \Box \]

\begin{prop}
If $\Phi_0(z)$ is the socle solution, then the function $\Psi(z)$ is
uniquely defined by the equations (\ref{exp-diff}) and
the condition $\Psi(0) = 0$ as an element of ${\bf C}[[z]]$.
\label{log-q}
\end{prop}

{\bf Proof. } Let $\Psi(z) = \sum_{i = - \infty }^{+ \infty} b_i z^i$
be an element of $z{\bf C}[[z ]]$, i.e., $b_i = 0$ for $i \leq 0$.
By \ref{solut}, for any  $n \in {\bf Z}$, the coefficient
by $z^n$ in ${\cal P}\Psi(z)$ is
\[ P_m(n)b_n + P_{m-1}(n-1)b_{n-1} + \cdots + P_0(n -m)b_{n-m}. \]
On the other hand, the coefficient by $z^n$ in ${\cal P}_{\Theta}'\Phi_0(z)$
is
\[ P_m'(n)a_n + P_{m-1}'(n-1)a_{n-1} + \cdots + P_0'(n -m)a_{n-m}. \]
Thus, we obtain the recurrent linear non-homogeneous relation
\begin{equation}
 P_m(n)b_n + \sum_{i =1}^{m} P_{m -i}(n-i) b_{n -i} +
\sum_{i =0}^{m} P_{m -i}'(n-i) a_{n -i} = 0.
\label{rec-psi}
\end{equation}

Since $P_m(n) = n^{d+1} \neq 0$ for $n \geq 1$, one can find all coefficients
$b_i$ $(i \geq 1)$ using (\ref{rec-psi}). For instance, we obtain
\[ b_1 = -(d+1)a_1 - P_{m-1}'(0)a_0 \]
\[ 2^{d+1}b_2 = -(d+1) 2^d a_2 - P_{m-1}(1)b_1 - P_{m-1}'(1)
a_1 - P_{m-2}'(0)a_0; \; ...\; {\rm etc.}  \]

\begin{coro}
Let ${\cal P}$ be a $MU$-operator then the quotient $\Psi/\Phi$ of
the solutions of the linear system
\[ {\cal P}\Phi =0,\;\; {\cal P}_{\Theta}' \Phi + {\cal P}\Psi,\;\;
\Psi(0) = 0 \]
is a function depending only on ${\cal P}$.
\end{coro}

We come now to the most important definition:

\begin{opr}
{\rm The element
\[ q = \exp \left( \frac{\Phi_1(z)}{\Phi_0(z)} \right)
= z \exp \left( \frac{\Psi(z)}{\Phi_0(z)} \right) \]
is called the {\em q-parameter} for the $MU$-operator ${\cal P}$.}
\end{opr}

\subsection{Generalized hypergeometric functions and $2$-term recurrent
relations}

\noindent

Since the number $m +1$ of terms in a recurrent relation (\ref{recurrent})
is at least $2$, $2$-term recurrencies are the simplest ones. Any such a
relations is defined by two polymomials $P_0(y)$ and $P_1(y)$ of degree
${d+1}$:

\begin{equation}
  P_0(n)a_{n} = P_1(n+1)a_{n +1}.
\label{rec-2}
\end{equation}

Without loss of generality we again assume that the leading coefficient
of $P_1(y)$ is $1$.
\begin{opr}
{\em Denote by
\[ G_{d+1}(\alpha,\beta;w) =
G_{d+1} \left( \begin{array}{c} \alpha_1, \ldots, \alpha_{d+1} \\
\beta_1, \ldots, \beta_{d+1} \end{array}; w  \right) \]
 the series
\[ \sum_{n \geq 0} \frac{ \prod_{i =1}^{d+1} \Gamma(\alpha_i)}
{\prod_{i =1}^{d+1} \Gamma(\beta_i)} \times
 \left(  \frac{  \prod_{i =1}^{d+1} \Gamma(\beta_i +n )}
{ \prod_{i =1}^{d+1} \Gamma(\alpha_i +n)} \right) w^n. \]
which is the generalized hypergeometric function
with parameters $\alpha_1,$ $\dots,$ $\alpha_{d+1}$, $\beta_1,$ $\ldots,$
$\beta_{d+1}$.
(This is a slight modification of the well-known generalized
hypergeometric function $\,_{d+1}F_{{d}}$ (see \cite{hyp.geom,slater}).)}
\end{opr}

\begin{prop}
Assume that
\[ P_1(y) = \prod_{i =1} (y + \alpha_i); \]
\[ P_0(y) = \lambda \prod_{i =1} (y + \beta_i) \]
then the function $G_{d+1}(\alpha,\beta; \lambda z)$ is a formal
solution of the differential equation
\[ {\cal P}\Phi = (P_1(\Theta) - z P_0(\Theta))\Phi = 0. \]
\end{prop}

Consider now the case when ${\cal D}$ is a $MU$-operator, i.e.,
$P_1(y) = y^{d+1}$, and the recurrent relation has the form
\[ (n+1)^{d+1} a_{n+1} = P_{0}(n) a_n. \]
Then for  the power series $\Psi(z) =
\sum_{i \geq 1} b_i z^i$ which is the solution to
\[ {\cal D}_{\Theta}' \Phi_0(z) +
{\cal P}  \Psi(z) = 0, \]
where \[ \Phi_0(z) = \sum_{i = 0}^{\infty} a_i z^i, \]
is a regular solution to ${\cal P}\Psi =0$,
the coefficients $\{ b_i \}$ satisfy the
recurrent relation
\[ n^{d+1}b_n = P_1(n-1)b_{n-1} + P_1'(n-1)a_{n-1} - (d+1)n^d a_n. \]

\subsection{$d$-point Yukawa functions}

\noindent

Let $\pi\; :\;  V_z  \rightarrow S$
be a $1$-parameter family of Calabi-Yau $d$-folds,  where
$S = {\rm Spec}\, {\bf C}[[z]]$. Let $T$ be the corresponding
monodromy transformation acting on $H_d(W_z, {\bf C})$, $T_u$ the
unipotent part part of $T$, $N = Log\, T_u$.

\begin{opr}
{\rm The family $V_z$ is said to have the {\em  maximal unipotent
monodromy at $z = 0$} if the weight filtration
\[ 0 = W_{-1} \subset W_0 \subset W_1 \subset \cdots \subset
W_{d-1} \subset W_d = H^d(V_z, {\bf C}) \]
defined by $N$ is orthogonal to the Hodge filtration $\{ F^i \}$, i.e.,
\[ H^d(V_z, {\bf C}) = W^{\perp}_{i} \oplus F^{d-i} \; i =0, \ldots , d.\]}
\end{opr}
(This is similar to definitions given in \cite{greene1,morrison.hodge}.)

Assume that the family $V_z / S$ has the maximal unipotent monodromy at
$z =0$ and ${\rm dim}\, F^i/F^{i+1} = 1$ for $ i = 0,\ldots, d$. Then the
Jordan normal form of $N$ has exactly {\em one} sell of size $(d+1)\times
(d+1)$. This means that there exists a $d$-cycle $\gamma \in H_d(V_z,
{\bf Z})$
such that $\gamma, N\gamma, \ldots, N^d \gamma$ are linearly independent in
$H_d(V_z, {\bf Z})$, and $N^d \gamma = \gamma_0$ is a monodromy invariant
$d$-cycle.
Take a $1$-parameter
family $\omega(z)$ of holomorphic $d$-forms on
$W_z$.
It is well-known that the periods of $\omega(z)$ over the
$d$-cycles in $H_d(V_z,{\bf C})$
satisfy a  Picard-Fuchs differential equation  of order $d+1$
defined by
some differential $MU$-operator

\begin{equation}
{\cal P} = \Theta^{d+1} + C_d(z) \Theta^d + \cdots + C_0(z).
\label{P-F}
\end{equation}

\begin{opr}
{\rm Define the coupling functions $W_{k,l}(z)$ $(k,l \geq 0, \; k,l \in {\bf
Z})$
as follows
\[ W_{k,l} = \int_{V_z} \Theta^k \omega(z) \wedge \Theta^l \omega(z). \]
(By definition, we put $\Theta^0 = 1$ to be the identity differential
operator.) }
\label{couplings}
\end{opr}

\begin{opr}
{\rm \cite{morrison.mirror} The coupling function $W_{d,0}$ is
called {\em unnormalized
$d$-point Yukawa function}.}
\end{opr}

\begin{prop}
The coupling functions $W_{k,l}(z)$ satisfy the properties

{\rm (i)} $W_{k,l}(z) = (-1)^d W_{l,k}$;

{\rm (ii)} $W_{k,l}(z) = 0$ for $k + l < d$;

{\rm (iii)}  $\Theta W_{k,l}(z) = W_{k+1,l}(z) + W_{k,l+1}(z)$.

{\rm (iv)} $W_{d+k +1,0}(z) + C_d(z) W_{d+k,0}(z) + \cdots +
C_0(z)W_{k,0} (z) = 0.$
\label{w-properties}
\end{prop}

{\bf Proof.} The statements follow immediatelly from the properties
of the cup-product and from the Griffiths
transversality property.

\begin{theo}
The $d$-pont Yukawa function  $W_{d,0}(z)$ satisfies the
linear differential equation
of order one
\begin{equation}
 \Theta W_{d,0}(z) + \frac{2}{d+1}C_d(z) W_{d,0} =0.
\label{coup-eq}
\end{equation}
\end{theo}

{\bf Proof. }  By \ref{w-properties}(ii), we have
\begin{equation}
 W_{d-i,i}(z) + W_{d-i-1,i+1}(z) =0\;\;{\rm  for}\; i =0,1, \dots, d.
\label{sum1}
\end{equation}
Therefore, $W_{d,0}(z) = (-1)^i W_{d-i,i}$. On the other hand,
by \ref{w-properties}(iii), we have
\begin{equation}
\Theta W_{d-i,i} = W_{d-i+1,i}(z) + W_{d-i,i+1}(z)
\;\;{\rm  for}\; i =0,1, \dots, d.
\label{sum2}
\end{equation}

It follows from (\ref{sum1}) and (\ref{sum2}) that
\begin{equation}
k \Theta W_{d,0}(z) = \sum_{i=0}^{k-1} (-1)^i W_{d-i,i}(z) =
 W_{d+1,0}(z) + (-1)^{k-1}  W_{d- k + 1,k}(z).
\label{sum3}
\end{equation}

{\sc Case I:} $d$ is odd. Since
\[ W_{\frac{d+1}{2}, \frac{d+1}{2}}(z) = 0\;\;
{\rm (\ref{w-properties}(i))},\]
we obtain
\[ \Theta W_{\frac{d+1}{2}, \frac{d-1}{2}}(z) =
W_{\frac{d+3}{2}, \frac{d-1}{2}}(z). \]
Using (\ref{sum1}) and (\ref{sum3}) for $k = (d-1)/2$, we obtain
\[ \frac{(d+1)}{2} \Theta W_{d,0}(z) =  W_{d+1,0}(z). \]
By \ref{w-properties}(ii) and (iv), this implies the equation
(\ref{coup-eq}) for $W_{d,0}(z)$.

{\sc Case I:} $d$ is odd. One has
\[ \Theta W_{\frac{d}{2}, \frac{d}{2}}(z) =
W_{\frac{d+2}{2}, \frac{d}{2}}(z) + W_{\frac{d}{2}, \frac{d+2}{2}}(z) =
2 W_{\frac{d+2}{2}, \frac{d}{2}}(z). \]

Using (\ref{sum1}) and (\ref{sum3}) for $k = d/2$, we obtain
\[ (d+1)\Theta W_{d,0}(z) = 2  W_{d+1,0}(z). \]
The latter again implies the same linear differential equation for
$W_{d,0}(z)$.

\begin{coro}
\[ W_{d,0}(z) = c_0\exp \left( {-\frac{2}{d+1} \int_0^z
C_d(v) \frac{dv}{v}}\right)\]
for some nonzero constant $c_0 = W_{d,0}(0)$.
\label{solution}
\end{coro}

\begin{exam}
{\rm Assume that ${\cal P} = \Theta^{d+1} - zP_0(\Theta)$ be the
$MU$-operator correspondng to a $2$-term recurrent relation
$(n+1)^{d+1} a_{n+1} = P_0(n)a_n$, where $P_0(y) = \lambda y^{d+1} + \cdots$.
Then the Yukawa $d$-point function $W_{d,0}(z)$ equals
\[ W_{d,0}(z) = \frac{c_0}{1 - \lambda z}, \]
i.e., $W_{d,0}(z)$ is a rational function in $z$. }
\label{yukawa}
\end{exam}

\subsection{Multidimensional Picard-Fuchs differential
systems with a symmetry group}

\noindent

So far we considered only the case of the $1$-parameter
family of Calabi-Yau $d$-folds $V_z$ such that
${\rm dim}\, F^i/F^{i+1} =1$ for $ i =0, \ldots, d$. It is easy to see
that the same methods can be applied to the case
${\rm dim}\, F^i/F^{i+1} \geq 1$ provided $V_z$ has a large authomorphisms
group.
\medskip

\begin{prop}
Let $V_z$ be a $1$-parameter family of Calabi-Yau $d$-folds with
${\rm dim}\, F^i/F^{i+1}$ $\geq 1$. Assume that there exists an
action of a finite group $G$ on $V_z$ such that
the $G$-invariant part $(F^i/F^{i+1})^G$ is $1$-dimensional for all
$i =0, \ldots, d$. Then
the holomorphic differential $d$-form $\omega(z)$ again satisfies
the Picard-Fuchs differential equation of order $d +1 $.
\label{symmetry}
\end{prop}

{\bf Proof. } The statement immediately follows from the fact that
the cohomology classes of $\omega(z), \Theta \omega(z), \ldots,
\Theta^d \omega(z)$ form the basis of the $G$-invariant subspace
$H^d(V_z, {\bf C})^G \subset H^d(V_z, {\bf C})$. \hfill $\Box$.
\bigskip

\section{Calabi-Yau complete intersections in ${\bf P}^N$}

\subsection{Rational curves and generalized hypergeometric series}

\noindent

Let $V$ be a complete intersection of $r$-hypersurfaces $V_1,
\ldots, V_r$ of degrees $d_1, \ldots, d_r$ in  ${\bf P}^{d+r}$.
Then $V$ is a Calabi-Yau $d$-fold if and only if $d + r + 1 =
d_1 + \ldots + d_r$. A rational curve
$C$ of degree $n$ in ${\bf P}^{d+r}$
has $nd_i$ intersection points with a generic hypersurface $V_i$.
On the other hand, there
exists a degeneration of every divisor $V_{i}$ into the union of
$d_i$ hyperplanes. Each of these hyperplanes has $n$ intersection
points with $C$.
This motivates the definition of the corresponding
generalized hypergeometric series  $\Phi_0(z)$ as

\begin{equation}
\sum_{i =0}^{\infty} \frac{(nd_1!)}{(n!)^{d_1}} \cdots \frac{(nd_r !)}
{(n!)^{d_r}} z^n
\label{proj-space}
\end{equation}

The coefficients
\[ a_n = \frac{(nd_1!)\cdots (nd_r !)}{ (n!)^{d+r +1}} \]
satisfy the recurrent relation
\[ (n+1)^{d+1} a_{n+1} = P(n) a_n \]
where $P(y)$ is the polynomial of degree $d +1$:
\[ P(n) = \frac{(nd_1 + d_1)!}{(nd_1) !}  \cdots
\frac{(nd_r + d_r)!}{(n d_r) !}(n+1)^{-r} = \lambda n^{d+1} + \cdots . \]
In particular, the leading coefficient of $P(y)$ is $\lambda =
\prod_{ i=1}^{r} (d_i)^{d_i}$.

\begin{exam}
{\rm Let $V$ be a complete intersection of two cubics in ${\bf P}^5$.
The corresponding generalized hypergeometric series is
\[ \Phi_0(z) = \sum_{n \geq 0} \frac{(3n!)^2}{(n!)^6} z^n. \]
This series was found in \cite{lib.teit} using the explicit construction
of mirrors for $V$ by orbifolding the $1$-parameter family of
special complete intersections of two cubics in ${\bf P}^5$:
\[ Y_1^3 + Y_2^3 + Y_3^3 = 3 \psi Y_4Y_5Y_6; \]
\[  Y_4^3 + Y_5^3 + Y_6^3 = 3 \psi Y_1Y_2Y_3, \]
by an abelian group $G$ of order $81$.
where $z = (3\psi)^{-6}$.
\medskip

We will give another interpretation of the construction of mirrors
$V'$ for $V$ which
immediatelly implies that $\Phi_0(z)$ is the monodromy invariant
period for the regular differential $3$-form on $V'$.

Let $Z_{f_1 f_2}$ be the complete intersection of two hypersurfaces in a
$5$-dimensional algebraic torus ${\bf T} = {\rm Spec} [ X_1^{\pm 1},
\ldots, X_5^{\pm 1} ]$ defined by the Laurent polynomials
\[ f_1(u,X) = 1 - (u_1 X_1 + u_2 X_2 + u_3 X_3), \;
f_2(u,X) = 1 - (u_4 X_4 + u_5X_5 + u_6(X_1 \cdots X_5)^{-1} ). \]
We define the differential $3$-form $\omega$ on
$Z_{f_1 f_2}$ as the residue of the rational differential $5$-form on
${\bf T}$:
\[ \omega = \frac{1}{(2 \pi \sqrt{-1})} {\rm Res}
\frac{ 1}{f_1(u,X)f_2(u,X)} \frac{dX_1}{X_1} \wedge \cdots \wedge
\frac{dX_5}{X_5}.\]
Let $z = u_1 \cdots u_6$. By residue theorem, we obtain
\[ \Phi_0(z) = \frac{1}{(2 \pi \sqrt{-1})} \int_{\mid X_i \mid =1}
\frac{ 1}{f_1(u,X)f_2(u,X)} \frac{dX_1}{X_1} \wedge \cdots \wedge
\frac{dX_5}{X_5}.\]
In this interpretation, the mirrors for $V$ are smooth Calabi-Yau
compactifications of of affine $3$-folds $Z_{f_1 f_2}$.

The equivalence between  the above two construction of mirrors for $V$ follows
by  the substitution
\[ X_1 = Y_1^3/(Y_3Y_4 Y_5), \;  X_2 = Y_2^3/(Y_3Y_4 Y_5), \;
 X_3 = Y_1^3/(Y_3Y_4 Y_5), \; \]
\[  X_4 = Y_4^3/(Y_1Y_2 Y_3), \;  X_5 = Y_5^3/(Y_1Y_2 Y_3), \;\]
\[ u_1 = \cdots = u_6 = (3\psi)^{-1}. \]  }
\label{3-3}
 \end{exam}

\begin{prop}
The normalized Yukawa $d$-differential for Calabi-Yau complete
intersections has the form
\[ {\cal W}_{d} = \frac{d_1 \cdots d_r}{ (1 - \lambda z)\Phi_0^2(z)}
\cdot (\frac{dz}{z})^{\otimes d},  \]
where $d_1, \ldots, d_i$ are degrees of hypersurfaces.
\end{prop}

{\bf Proof. } The statement follows from
\ref{yukawa} and the normalizing condition  $d_1 \cdots d_r = W_d(0)$.
\hfill $\Box$

\subsection{The construction of mirrors}

\noindent

Let $V$ be a $d$-dimensional
Calabi-Yau complete intersection of $r$ hypersurfaces
of degrees $d_1, \ldots, d_r$ in ${\bf P}^{d + r}$. We propose
the explicit construction of candidats for mirrors with respect to
$V$ as follows:

Let $E = \{ v_1, \ldots, v_{d +r +1}\} $ be a generating set in
the $(d+r)$-dimensional lattice $N \cong {\bf Z}^{d +r}$ such that there exist
the relation
\[ v_1 +  \cdots +  v_{d +r + 1} =0. \]
We divide $E$ into a disjoint union of $r$ subsets $E_i \subset E$ such that
${\rm Card}\, E_i = d_i$.
For $i =1, \ldots, r$, we define the Laurent polynomial
$P_i(u, X)$ in variables $X_1, \ldots, X_{d+r}$ as
\[ P_i(X) = 1 - (\sum_{v_j \in E_i} u_j X^{v_j}), \]
where $u_1, \ldots, u_{d +r +1}$ are independent parameters.
We denote by  $V'$  a  Calabi-Yau compactification of
$d$-dimensional affine complete intersections $Z$ in ${\bf T} =
{\rm Spec}\, [ X_1^{\pm 1}, \ldots, X_{d +r}^{\pm1}]$  defined
by the polynomials $P_1(u,X), \ldots, P_r(u,X)$ with sufficiently
general coefficients $u_i$. It is easy to see that up to an isomorphism
the affine variaties $Z \subset {\bf T}$ depend only on
$z = u_1 \cdots u_{d+r+1}$. Thus, we have obtained a $1$-parameter family
of $d$-dimensional varieties $V'$.

\begin{conj}
The $1$-parameter family of $d$-dimensional varieties $V'$ yields the mirror
family for $V$.
\end{conj}

This conjecture is motivated by the combinatorial interpretation proposed
in \cite{batyrev.dual} of the well-known construction of mirrors
for hypersurfaces of degree $d+2$ in ${\bf P}^{d+1}$ (see \cite{greene}).
On the other hand, the conjecture is supported by the following property:

\begin{prop}
The hypergeometric series $\Phi_0(z)$ in {\rm (\ref{proj-space})}
is the monodromy invariant period function of the holomorphic $d$-form
$\omega$ on $V'$.
\end{prop}

{\bf Proof. } The statement follows from the equality
\[\sum_{i =0}^{\infty} \frac{(nd_1!)}{(n!)^{d_1}} \cdots \frac{(nd_r !)}
{(n!)^{d_r}} z^n  =  \frac{1}{(2\pi \sqrt{-1})^{d +r}}
\int_{\mid X_j \mid =1}  \frac{1}{P_{1}(X) \cdots P_{r}(X)}
\frac{dX_1}{X_1} \wedge \cdots \wedge \frac{dX_{d+r}}{X_{d+r}}. \]
\hfill $\Box$.
\bigskip

\section{Complete intersections in toric varieties}

\subsection{The generalized hypergeometric series $\Phi_0$}

\noindent

Let $N$ be a free abelian group of rank $(d+r)$. Consider $r$  finite
sets
\[ E_i = \{ v_{i, 1}, \ldots, v_{i,k_i} \},\; i =1, \ldots, r \]
consisting of elements $v_{i,j} \in N$. Let $E$ be the union
$E_1 \cup \cdots \cup E_r$.

We put $ k = {\rm Card}\, E = k_1 + \ldots + k_r$ and assume that $E$
generates the group $N$.
Let $R(E)$ be the subgroup in ${\bf Z}^n$ consisting of
all integral vectors $\lambda = \{ \lambda_{i,j} \}$ such that
\[ \sum_{i=1}^r \sum_{j =1}^{k_i} \lambda_{i,j} v_{i,j} = 0. \]
We denote by $R^+(E)$ submonoid in $R(E)$ consisiting of all
$\lambda = \{ \lambda_{i,j} \} \in R(E)$ such that $\lambda_{i,j} \geq 0$.

\begin{opr}
{\rm Let $u_{i,j}$ be $k$ independent complex variables parametrized
by $k$ integral vectors $v_{i,j}$. Define the power series $\Phi_0(u)$ as

\[ \Phi_0(u) = \sum_{\lambda \in R^+(E)}
 \prod_{i =1}^r (\sum_{j =1}^{k_i}
\lambda_{i,j})!  \left( \prod_{j =1}^{k_i}
 \frac{u_{i,j}^{\lambda_{i,j}}}{(\lambda_{i,j})!}\right) .\]}
\end{opr}

Let $\lambda^{(1)}, \ldots, \lambda^{(t)}$ be a ${\bf Z}$-basis of the lattice
$R(E)$ such that every element $\lambda \in R^+(E)$ is
a non-negative integral linear combination of $\lambda^{(i)}$.
We define new $r$ complex variables $z_1, \ldots, z_s$ as follows
\[ z_s =  \prod_{i =1}^t \prod_{j =1}^{k_i} u_{i,j}^{\lambda_{i,j}^{(s)}};
\; s =1 , \ldots, t. \]
Thus, the series $\Phi_0(u)$ can be rewritten as the power series
$\Phi_0(z)$ in $t$ variables $z_1, \ldots, z_t$.
\medskip

\begin{exam}
{\rm Let $E = \{ v_1, \ldots, v_{d+1} \}$ be vectors generating
$d$-dimensional lattice $N$ and satisfying the integral relation
$v_1 + \cdots + v_{d+1} = 0$, i.e., the group $R(E)$ is generated
by the vector $(1, \ldots, 1) \in {\bf Z}^{d+1}$. Then the corresponding
generalized hypergeometric series is
\[ \Phi_0(u) = \sum_{n \geq 0} \frac{(nd +n)!}{(n!)^{d+1}}
(u_1 \cdots u_{d+1})^n =
\sum_{n \geq 0} \frac{(nd +n)!}{(n!)^{d+1}}
z^n = \Phi_0(z), \]
where $z = u_1 \cdots u_{d+1}$.  The integral representation of this
series is the monodromy invariant period function for mirrors of
hypersurfaces of degree $(d+1)$ in ${\bf P}^d$.
}
\end{exam}

\begin{opr}
{\rm Let ${\bf T}$ be a $(d +r)$-dimensional algebraic torus with
the Laurent coordinates $X = (X_1, \ldots, X_{d+r})$. We define
$r$ Laurent polynomials $P_{E_1}(X), \ldots, P_{E_r}(X)$ as follows
\[ P_{E_i}(X) = 1 - \sum_{v_{i,j} \in E_i} u_{i,j} X^{v_{i,j}}. \]  }
\end{opr}

\begin{prop}
The series $\Phi_0(u)$ admits the following integral representation
\[ \Phi_0(u) = \frac{1}{(2\pi \sqrt{-1})^{d +r}}
\int_{\mid X_j \mid =1}  \frac{1}{P_{E_1}(X) \cdots P_{E_r}(X)}
\frac{dX_1}{X_1} \wedge \cdots \wedge \frac{dX_{d+r}}{X_{d+r}}. \]
\label{integral}
\end{prop}

{\bf Proof. } The statement follows immediately from the residue formula.
\hfill $\Box$

\subsection{Calabi-Yau complete intersections}

\noindent

Let ${\bf P}_{\Sigma}$ be a quasi-smooth $(d+r)$-dimensional projective
toric
variety defined
by a $(d +r)$-dimensional simplicial fan $\Sigma$.
Assume that there exist $r$
line bundles
${\cal L}_1, \ldots, {\cal L}_r$ such that each ${\cal L}_i$ is
generated by global sections and the tensor product
\[ {\cal L}_1 \otimes  \cdots \otimes  {\cal L}_r \]
is isomorphic to the anticanonical bundle on ${\cal K}^{-1}$ on
${\bf P}_{\Sigma}$.
If $V_i$ is the set of zeros  of a generic global section of ${\cal L}_i$,
then the complete intersection
$V = V_1 \cap \cdots \cap V_r$ is a $d$-dimensional Calabi-Yau variety
having only Gorenstein toroidal singularities which are analytically
isomorphic to toric singularities of ${\bf P}_{\Sigma}$.
\medskip

Now let $E= \{v_1, \ldots, v_k\}$ be the set of all generators of
$1$-dimensional cones in
$\Sigma$. Denote by $D_j$ the toric divisor on ${\bf P}_{\Sigma}$
corresponding to $v_j$.
Notice that
\[ {\cal K}^{-1} = \otimes_{j=1}^k {\cal O}_{{\bf P}_{\Sigma}}(D_j). \]
Following a suggestion of Yu. I. Manin, we assume that one can represent $E$ as
a disjoint union
  \[ E = E_1 \cup \cdots \cup E_r \]
such that the line bundle ${\cal L}_i$ is isomorphic to
the tensor product
\[ \otimes_{v_j \in E_i} {\cal O}_{{\bf P}_{\Sigma}}(D_j). \]

The elements of the group $R(E)$ can be identified with the homology classes
of $1$-dimensional algebraic cycles on ${\bf P}_{\Sigma}$.
Moreover, one has the following property

\begin{prop}
Let $\lambda = (\lambda_1, \ldots, \lambda_k)$  be an arbitrary element of
$R(E)$ representing the class of an algebraic $1$-cycle $C$. Then
\[ \lambda_i = \langle D_i, C \rangle,  \;i =1, \ldots, k. \]
\end{prop}

We can always choose
a ${\bf Z}$-basis $\lambda^{(1)}, \ldots, \lambda^{(t)}$ of $R(E)$
such that every {\em effective} $1$-cycle on ${\bf P}_{\Sigma}$ is a
non-negative linear combination of the
elements  $\lambda^{(1)}, \ldots, \lambda^{(t)}$.
Since the submonoid $R^+(E)$ consists of classes of {\em nef-curves},
this implies that every element of $R^+(E)$ is also a non-negative
linear combination of $\lambda^{(1)}, \ldots, \lambda^{(t)}$.
This allows us to rewrite the series $\Phi_0(u)$ in $t$ algebraically
independent variables $z_1\, \ldots, z_t$ $( t = {\rm rk}\, R(E))$.

\begin{coro}
The series $\Phi_0(z)$ can be interpreted via the intersection numbers
of classes $\lbrack C \rbrack$ of
curves $C$ on ${\bf P}_{\Sigma}$ as follows
\[ \Phi_0(z) = \sum_{\lbrack C \rbrack \in R^+(E)}
\frac{ (\langle V_1,C \rangle)! \cdots (\langle V_r, C \rangle)!}
{\langle D_1, C \rangle ! \cdots \langle D_k, C \rangle !} z^{\lbrack C
\rbrack}, \]
where $z^{\lbrack C \rbrack} = z_1^{c_1} \ldots z_t^{c_t}$,
$ \lbrack C \rbrack = c_1 \lambda^{(1)} + \cdots + c_t \lambda^{(t)}$.
\end{coro}

\subsection{General conjectures}

\noindent

Let $V$ be a $d$-dimensional Calabi-Yau complete intersection
of hypersurfaces $V_1, \ldots, V_r$ in a $(d+r)$-dimensional
quasi-smooth toric variety defined by a simplicial fan $\Sigma$.
Choose a ${\bf Z}$-basis $\lambda^{(1)}, \ldots, \lambda^{(t)}$
in $R(E)$  such that the classes of all effective $1$-cycles have
non-negative integral coordinates.
We assume that the divisors $V_1, \ldots, V_r$ are numerically effective
(in particular, they are not assumed to be necessary ample divisors).
We assume also that the following conditions are satisfied:

{\rm (i)} $V$ is smooth;

{\rm (ii)} the restriction mapping ${\rm Pic}\, {\bf P}_{\Sigma} \rightarrow
{\rm Pic}\, V$ is injective.
\medskip

In this situation, there exist  two  flat A-modle connections:
the connection
$\nabla_{AP}$ on $H^*({\bf P}_{\Sigma})$ and the  connection
$\nabla_{AV}$ on
$H^*(V, {\bf C})$. Let $\tilde{H}^i$ be the image of $H^i({\bf P}_{\Sigma},
{\bf C})$ in $H^i(V, {\bf C})$.
The connection $\nabla_{AP}$ defines the quantum variation on
cohomology of toric variety ${\bf P}_{\Sigma}$. It follows from the result in
\cite{batyrev.quant} the following.

\begin{prop}
The complex coordinates $z_1, \ldots, z_t$ on $\tilde{H}^2$ can be identified
with flat coordinates with respect to $\nabla_{AP}$.
\end{prop}

\begin{conj}
The  generalized hypergeometric series
$\Phi_0(z)$ as a function of $\nabla_{AP}$-flat $z$-
coordinates
on $\tilde{H}^2$ is a solution of the differential system ${\cal D}$
defined by the A-model connection $\nabla_{AV}$ on $\tilde{H}^2$ which
defines the  quantum  variation of Hodge
structures on
the subring in $\bigoplus_{i = 0}^{d} H^{2i}(V, {\bf C})$ generated by
restrictions of the classes in ${\rm Pic}\, {\bf P}_{\Sigma}$ to
$V$.
\label{conj.1}
\end{conj}

\begin{rem}
{\rm One can check in many  examples that the differential system ${\cal D}$
is already  defined  by the generalized hypergeometric series
$\Phi_0(z)$. Probably there exists a general explanation of this fact.}
\end{rem}

\begin{conj}
The $\nabla_{AV}$-flat coordinates $q_1, \ldots, q_t$ on $\tilde{H}^2$
are
defined as
\[ q_i = \exp (\Phi_i(z)/\Phi_0(z)), \; i = 1, \ldots, t,  \]
where $\Phi_i(z)$ is a logarithmic
solution to the differential system ${\cal D}$ having
the form
\[ \Phi_i(z) = (\log z_i) \Phi_0(z) + \Psi_i(z), \; \Psi_i(0) = 0  \]
for some regular at $z =0$ power series $\Psi_i(z)$.

Moreover, all coefficients of the  expansion of $\nabla_{AV}$-flat
coordinates  $q_i$ as power series of $\nabla_{AP}$-flat $z$-coordinates are
{\bf integers}.
\label{conj.2}
\end{conj}

\begin{rem}
{\rm This conjecture establishes a general method for normalizing
the logarithimic solutions defining the canonical
$q$-coordinates for the differential system ${\cal D}$. There are two
motivations for this conjecture. First, the conjecture is true for all
already known examples of $q$-coordinates for Picard-Fuchs equations
corresponding to Calabi-Yau complete intersections in products of
weighted projective spaces (see examples in the remaining part of the
paper). Second, the Lefschetz theorem and the calculation of the quantum
cohomolgy ring of toric varieties \cite{batyrev.quant} imply  the relation
\[ q_i = z_i + O(\mid z \mid^2),\; ( i =1, \ldots, t). \]   }
\end{rem}
\medskip

\begin{conj}
Assume that $V$ has dimension $3$. Let $K_{i,j,k}(z)$ be structure
constants defining the A-model $\nabla_{AV}$
connection in the $z$-coordinates. Then
\[ \Phi_0^2(z)  K_{i,j,k}(z) \]
is a rational function in $z$-coordinates.
\label{conj.21}
\end{conj}

\begin{conj}
The mirror Calabi-Yau varieties with respect to $V$ are Calabi-Yau
compactifications of the complete intersection of the affine
hypersurfaces in the $(d+r)$-di\-men\-sio\- nal algebraic torus ${\bf T}$
defined
by the equations
\[ P_{E_1}(X) = \cdots = P_{E_r}(X) = 0.\]
\label{conj.3}
\end{conj}

\begin{rem}
{\rm Recall that two Calabi-Yau $d$-folds $V$ and $V'$ are called
mirror symmetric if $h^{p,d-p}(V) = h^{d-p,d-p}(V')$ and the superconformal
field theories corresponding to $V$ and $V'$ are isomorphic. It was
proposed in \cite{batyrev.dual} a general method for constructing pairs
of mirror symmetric Calabi-Yau hypersurfaces in toric varieties
based on the duality among so called {\em reflexive polyhedra} $\Delta$ and
$\Delta^*$. However,
the equality $h^{1,1}(\hat{Z}_f) = h^{d-1,1}(\hat{Z}_g)$ for the pair
of Calabi-Yau $d$-folds $\hat{Z}_f$ and $\hat{Z}_g$ corresponding to
the polyhedra $\Delta$ and $\Delta^*$ are not sufficient to prove the
mirror duality between $\hat{Z}_f$ and $\hat{Z}_g$ in full strength.
One needs to prove more: the isomorphism between the quantum cohomology of
$\hat{Z}_f$ and $\hat{Z}_g$. Since the quantum cohomology are defined by
the canonical form of the A-model connection $\nabla_A$ in $q$-coordinates,
Conjecture \ref{conj.1} and Proposition \ref{integral} yield more
evidence for validity of Conjecture \ref{conj.3}.
We give below one example showing
that Conjecture \ref{conj.3} agrees with an orbifold construction of mirrors
for complete intersections in product of projective spaces inspired
by superconformal field theories. }
\end{rem}

\begin{exam}
{\rm Let $V$ be a Calabi-Yau complete intersection of two hypersurfaces
of degrees $(3,0)$ and $(1,3)$ in the product ${\bf P}^3 \times {\bf P}^2$.

It is known that the mirrors for $V$ can be obtained by orbifolding
the complete intersection of two special hypersurfaces
\[ S_1 T_1^3 + S_2 T_2^3 + S_3 T_3^3  = \phi S_4 T_1 T_2 T_3, \]
\[ S_1^3 + S_2^3 + S_3^3 + S_4^3 = \psi S_1 S_2 S_3 \]
by the group
$G$ of order $27$,
where $(S_1: S_2 : S_3 :S_4)$ and $(T_1 : T_2 : T_3)$ are the homogeneous
coordinates on ${\bf P}^3$ and ${\bf P}^2$ respectively.

On the other hand, the $5$-dimensional
fan $\Sigma$ defining ${\bf P}^3 \times {\bf P}^2$ has
$7$ generators $\{ v_1,$ $\ldots,$ $v_7 \}$ $= E$ satisfying the relations
\[ v_1 + v_2 + v_3 + v_6 = v_4 + v_5 + v_7 = 0. \]
We choose vectors $v_1, \ldots, v_5$ as the basis of the $5$-dimensional
lattice $N$. The complete intersection $V$ is defined by dividing $E$ into
two subsets $E_1 = \{v_1, v_2, v_3 \}$ and $E_2 =\{ v_4, v_5, v_6, v_7 \}$.
The corresponding polynomials $P_{E_1}(X)$ and $P_{E_2}(X)$ are
\[ P_{E_1} = 1 - (u_1 X_1 + u_2 X_2 + u_3 X_3), \]
\[ P_{E_2} = 1 - (u_4 X_4 + u_5 X_5 + u_6(X_1 X_2 X_3)^{-1} +
u_7(X_4 X_5)^{-1}). \]

We obtain the equivalence between two construction of mirrors by
putting $u_1 = u_2 = u_3 = \phi^{-1}$, $u_4 = u_5 = u_6 = u_7 = \psi^{-1}$,
and
\[ X_1 = \frac{S_1 T_1^2}{S_4  T_2 T_3}, \;
X_2 = \frac{S_2 T_2^2}{S_4 T_1  T_3}, \;
X_3 = \frac{S_3T_3^2}{S_4 T_1 T_2 }, \]
\[ X_4 = \frac{S_1^2}{S_2S_3},\; X_5 = \frac{S_2^2}{S_1S_3}.  \]}
\end{exam}
\bigskip

\subsection{Calabi-Yau $3$-folds with $h^{1,1} = 1$}

\noindent

We consider below examples of the generalized
hypergeometric series corresponding to smooth
Calabi-Yau  complete intersections $V$ of $r$ hypersurfaces
in a toric variety ${\bf P}_{\Sigma}$ such that $h^{1,1}(V) = 1$.
By Lefschetz theorem,
$h^{1,1}({\bf P}_{\Sigma})$ must be also $1$. So  $\Sigma$ is
a $(r +3)$-dimensional fan with $(r +4)$ generators. There exists the unique
primitive integral linear relation
$\sum \lambda_i v_i = 0$ among the generators $\{v_i \}$ of $\Sigma$, i.e.,
${\rm rk}\, R(E) =1$ where $E = \{ v_i \}$ is the whole set of generators
of $\Sigma$
(${\rm Card}\, E = r + 4$).
\medskip

In all these examples the $MU$-operator ${\cal P}$ has form
\[ {\cal P} = \Theta^4 - \mu z (\Theta + \alpha_1)(\Theta + \alpha_2)
(\Theta + \alpha_3)(\Theta + \alpha_4), \]
where the numbers $\alpha_1, \ldots, \alpha_4$ are positive rationals
satisfying  the relations
\[ \alpha_1 + \alpha_4 = \alpha_2 + \alpha_3 = 1. \]
The Yukawa $3$-differential in $z$-coordinate has form
\[ {\cal W}_3 = \frac{W(0)}{(1 - \mu z)\Phi_0^2(z)}
\left(\frac{dz}{z}\right)^{\otimes 3}. \]

\begin{exam} {\sl Hypersurfaces in weighted projective spaces:} {\rm
In this case we obtain Ca\-la\-bi-Yau hypersurfaces
in the following weighted projective spaces ${\bf P}(\lambda_1, \cdots,
\lambda_5)$

\begin{center}

\begin{tabular}{|c|c|c|c|c|}  \hline
$(\lambda_1, \ldots, \lambda_5)$ & $\Phi_0(z)$ & $ W(0)$ & $\mu$ &
$(\alpha_1, \alpha_2,\alpha_3, \alpha_4)$ \\ \hline
$(1,1,1,1,1)$ &
${\displaystyle \sum_{n \geq 0} \frac{(5n)!}{(n!)^5}z^n}$ &
5 & $5^5$ & $(1/5,2/5,3/5,4/5)$  \\ \hline
$(2,1,1,1,1)$ &
${\displaystyle \sum_{n \geq 0} \frac{(6n)!}{(n!)^4(2n!)}z^n}$ &
$ 3 $ & $2^5 3^6$ & $(1/6,2/6,4/6,5/6)$  \\ \hline
$(4,1,1,1,1)$ & ${\displaystyle
\sum_{n \geq 0} \frac{(8n)!}{(n!)^4(4n!)}z^n}$ &
$ 2 $ & $2^{18}$ & $(1/8,3/8,5/8,7/8)$   \\ \hline
$(5,2,1,1,1)$ & ${\displaystyle
\sum_{n \geq 0} \frac{(10n)!}{(n!)^3(2n!)(5n!)}z^n}$ &
$ 1 $& $2^9 5^6$ & $(1/10,3/10,7/10,9/10)$  \\ \hline
\end{tabular}

\end{center}

The $q$-expansion of the Yukawa $3$-point function and  predictions
$n_d$ for number of rational curves on these hypersurfaces were obtained
in \cite{morrison.picard,klemm1,font}.}
\end{exam}
\medskip

\begin{exam} {\sl  Complete intersections in ordinary projective spaces:}
{\rm Let $V_{d_1, \ldots, d_r}$ denotes the complete intersection of
hypersurfaces
of degrees $d_1, \ldots, d_r$.

\begin{center}

\begin{tabular}{|c|c|c|c|c|}  \hline
   & $\Phi_0(z)$ & $W(0)$ & $\mu$ &
   $(\alpha_1, \alpha_2,\alpha_3, \alpha_4)$ \\ \hline
$V_{3,3} \subset {\bf P}^5$ &
${\displaystyle \sum_{n \geq 0} \frac{((3n)!)^2}{(n!)^6}z^n}$ &
$9$ & $3^6$ & $(1/3,1/3,2/3,2/3)$  \\ \hline
$V_{2,4} \subset {\bf P}^5 $ & ${\displaystyle
\sum_{n \geq 0} \frac{(2n)!(4n)!}{(n!)^6}z^n}$ &
$8$ & $2^{10}$ & $(1/4,2/4,2/4,3/4)$  \\ \hline
$V_{2,2,3} \subset {\bf P}^6$ & ${\displaystyle
\sum_{n \geq 0} \frac{((2n)!)^2(3n)!}{(n!)^7}z^n}$ &
$12$ & $2^4 3^3$ & $(1/3,1/2,1/2,2/3)$  \\ \hline
$V_{2,2,2,2} \subset {\bf P}^8$ & ${\displaystyle
\sum_{n \geq 0} \frac{((2n)!)^4}{(n!)^8}z^n}$ &
$16$ & $2^8$ & $(1/4,1/4,1/4,1/4)$  \\ \hline
\end{tabular}

\end{center}

These  Calabi-Yau complete intersections in ordinary projective
spaces were  considered by Libgober and Teitelbaum  \cite{lib.teit}. }
\end{exam}
\medskip

\begin{exam} {\sl  Complete intersections in weighted projective spaces:}
{\rm
\begin{center}

\begin{tabular}{|c|c|c|c|c|}  \hline
   & $\Phi_0(z)$ & $W(0)$ & $\mu$ &
   $(\alpha_1, \alpha_2,\alpha_3, \alpha_4)$ \\ \hline
$V_{4,4} \in {\bf P}(1,1,1,1,2,2)$ &  ${\displaystyle \sum_{n =0}^{\infty}
\frac{(4n!)^2}{(n!)^4 (2n!)^2} z^n  }$ & $4$ & $2^{12}$ &
$(1/4,1/4,3/4,3/4)$ \\ \hline
$V_{6,6} \in {\bf P}(1,1,2,2,3,3)$ &  ${\displaystyle \sum_{n =0}^{\infty}
\frac{(6n!)^2}{(n!)^2 (2n!)^2(3n!)^2}z^n }$ & $1$ & $2^8 3^6$ &
$(1/6,1/6,5/6,5/6)$ \\ \hline
$V_{3,4} \in {\bf P}(1,1,1,1,1,2)$ & ${\displaystyle \sum_{n =0}^{\infty}
\frac{(4n!)(3n!)}{(n!)^5 (2n!)} z^n}$ &
$ 6$ & $2^6 3^3$ & $(1/4,1/3,2/3,3/4)$  \\ \hline
$V_{2,6} \in {\bf P}(1,1,1,1,1,3)$ & ${\displaystyle \sum_{n =0}^{\infty}
\frac{(6n!)(2n!)}{(n!)^5 (3n!)} z^n }$ & $4$ & $2^8 3^3$ &
$(1/6,1/2,1/2,5/6)$ \\ \hline
$V_{4,6} \in {\bf P}(1,1,1,2,2,3)$ & ${\displaystyle \sum_{n =0}^{\infty}
\frac{(6n!)(4n!)}{(n!)^3 (2n!)^2(3n!)} z^n}$ & $2$ & $2^{10}3^3$ &
$(1/6,1/4,3/4,5/6)$ \\ \hline
\end{tabular}

\end{center}

The coefficients of the Yukawa 3-point function $K_q^{(3)}$
for these five examples of Calabi-Yau $3$-folds $V$ having the  Hodge
number $h^{1,1}(V) =1$  were obtained by
A. Klemm and S. Theisen \cite{klemm2}. }
\end{exam}
\bigskip

\section{Calabi-Yau 3-folds in ${\bf P}^2 \times {\bf P}^2$}

\subsection{The generalized hypergeometric series $\Phi_0$}

\noindent

Calabi-Yau $3$-folds $V$ in ${\bf P}^2 \times {\bf P}^2$ are hypersurfaces of
degree $(3,3)$. The homology classes of rational curves on
${\bf P}^2 \times {\bf P}^2$ are parametrized by pairs of integers
$(l_1, l_2)$. Let $\gamma_1$, $\gamma_2$ be the homology classes of
$(1,0)$-curves and $(0,1)$-curves respectively. Then for any K\"ahler
class $\eta$ we put
\[ z_i = \exp (- \int_{\gamma_i} \eta ),\;\;(i =1,2) . \]

The generalized hypergeometric series corresponding to the fan $\Sigma$
defining ${\bf P}^2 \times {\bf P}^2$ is

\[ \Phi_0(z_1,z_2) = \sum_{l_1,l_2 \geq 0}
\frac{(3l_1 + 3l_2)!}{(l_1!)^3 (l_2!)^3} z_1^{l_1}z_2^{l_2}. \]

There are obvious two recurrent relations for the coefficients
$a_{l_1,l_2}$ of the series
\[ \Phi_0(z_1,z_2)  = \sum_{l_1,l_2 \geq 0}
a_{l_1,l_2} z_1^{l_1}z_2^{l_2} :\]

\[ (l_1 + 1)^3 a_{l_1+1,l_2} = (3l_1 + 3l_2 +1)(3l_1 + 3l_2 +2)
(3l_1 + 3l_2 +3)a_{l_1, l_2}; \]
\[ (l_2 + 1)^3 a_{l_1,l_2+1} = (3l_1 + 3l_2 +1)(3l_1 + 3l_2 +2)
(3l_1 + 3l_2 +3)a_{l_1, l_2}; \]

 Let
\[ \Theta_1 = z_1  \frac{\partial}{\partial z_1}, \;\;
\Theta_2 = z_2 \frac{\partial}{\partial z_2} \]
Then the  function $\Phi_0(z_1,z_2)$
satisfies the Picard-Fuchs differential system
${\cal D}$:
\[ \left( \Theta_1^3 - z_1 (3 \Theta_1 + 3\Theta_2 + 1 )
( 3\Theta_1 + 3\Theta_2 + 2 )( 3\Theta_1 + 3\Theta_2 + 3 ) \right)
\Phi_0 = 0 ; \]
\[ \left( \Theta_2^3 - z_2 ( 3 \Theta_1 + 3 \Theta_2 + 1 )
( 3 \Theta_1 + 3\Theta_2 + 2 )( 3 \Theta_1 + 3\Theta_2 + 3 ) \right)
\Phi_0 = 0. \]

The differential system ${\cal D}$ has the maximal unipotent monodromy
at $(z_1,z_2) =(0,0)$. There are two uniquely determined regular
at $(0,0)$ functions $\Psi_1(z_1,z_2)$ and
$\Psi_2(z_1,z_2)$ such that
\[ (\log z_1) \Phi_0(z_1,z_2) + \Psi_1(z_1,z_2), \]
\[ (\log z_2) \Phi_0(z_1,z_2) + \Psi_2(z_1,z_2) \]
are solutions to ${\cal D}$, and $\Psi_1(0,0) = \Psi_2(0,0) =0$.
If we put
\[ \Psi_j(z_1,z_2) = \sum_{ \begin{array}{c} {\scriptstyle l_1,l_2 \geq 0}\\
{\scriptstyle (l_1, l_2) \neq (0,0)} \end{array}}
b_{l_1,l_2}^{(j)} z_1^{l_1} z_2^{l_2}, \]
then one finds the coefficients $b_{l_1,l_2}^{(j)}$ from the simple
recurrent relations based on \ref{solut}.

The $q$-coordinates $q_1$, $q_2$ defined by the formulas
\[ q_1 = z_1 \exp (\Psi_1/\Phi_0), \;q_2 = z_2 \exp (\Psi_2/\Phi_0) \]
are the power series with integral coefficients in $z_1, z_2$ of the form
\[ q_j (z_1,z_2) = z_j \left( 1 + \sum_{ \begin{array}{c}
{\scriptstyle l_1,l_2 \geq 0}\\
{\scriptstyle (l_1, l_2) \neq (0,0)} \end{array}}
c_{l_1,l_2}^{(j)} z_1^{l_1} z_2^{l_2} \right), \; j =1,2. \]
By symmetry, one has $c_{l_1,l_2}^{(1)} = c_{l_2,l_1}^{(2)}$.

\subsection{Mirrors and the discriminant}

\noindent

Let $f$  be the Laurent polynomial
\[ f(X,u) = 1 - u_1X_1 - u_2X_2 - u_3(X_1X_2)^{-1} - u_4 X_3 - u_5X_4
-u_6 (X_3 X_4)^{-1}. \]
Let $\gamma_0$ be a generator of $H_4(({\bf C}^*)^4, {\bf Z})$, i.e.,
the cycle defined by the condition $\mid X_i \mid =1$ for $i =1,\ldots 4$.

By the residue theorem, the integral
\[ I(u) = \frac{1}{(2\pi \sqrt{-1})^4}
\int_{\gamma_0}  \frac{1}{f(X)}\frac{dX_1}{X_1}\wedge
\frac{dX_2}{X_2}\wedge \frac{dX_3}{X_3}\wedge \frac{dX_4}{X_4}\]
is the power series
\[ I(u) =  \sum_{k,m \geq 0}
\frac{(3k + 3 m)!}{(k!)^3 (m!)^3} (u_1 u_2 u_3)^k
(u_4 u_5 u_6)^m. \]
Thus, putting $z_1 = u_1 u_2 u_3$; $z_2 = u_4 u_5 u_6$, we obtain exactly
the generalized hypergeometric function $\Phi_0(z_1,z_2)$.

It was proved in \cite{batyrev.dual} and \cite{batyrev.var} that the function
$I(u)$ can be considered as the monodromy invariant period of the holomorphic
differential $3$-form
\[ \omega = \frac{1}{(2\pi \sqrt{-1})^4}
{\rm Res}\; \frac{1}{f(X)}\frac{dX_1}{X_1}\wedge
\frac{dX_2}{X_2}\wedge \frac{dX_3}{X_3}\wedge \frac{dX_4}{X_4}.\]
for the family of Calabi-Yau 3-folds $\hat{Z}_f$ which  are smooth
compactifications of the affine hypersurfaces $Z_f$ in $({\bf C}^*)^4$ defined
by Laurent polynomial $f$. One has $h^{1,1}(\hat{Z}_f)= 83$,
$h^{2,1}(\hat{Z}_f) = 2$.
The coordinates $z_1, z_2$ are natural coordinates on the moduli space of
Calabi-Yau $3$-folds $\hat{Z}_f$.
\bigskip

The mirror construction helps to understand the  the discriminant of
the differential system ${\cal D}$  as a polynomial function
in $z_1, z_2$.

By definition  \cite{gelfand}, the zeros of the discriminat are exactly
those values of the coefficients $\{u_i \}$ of $f(X)$ such that
the system
\[ f(X) = X_1 \frac{\partial }{\partial X_1}f(X) =
X_2 \frac{\partial }{\partial X_2}(X) =
X_3 \frac{\partial }{\partial X_3}f(X) = X_4
\frac{\partial }{\partial X_4}f(X) = 0 \]
has a solution in the toric variety ${\bf P}_{\Delta}$, where $\Delta$ is
the Newton polyhedron of $f$. Since ${\bf P}_{\Delta}$ is
isomorphic to the subvariety of ${\bf P}^6$ defined as
\[ {\bf P}_{\Delta} = \{ ( Y_0 : \ldots : Y_6) \in {\bf P}^6 \mid
Y_0^3 = Y_1 Y_2 Y_3, Y_0^3 = Y_4 Y_5 Y_6 \}, \]
or equivalently,
the system of the homogeneous equations
\[  u_0 Y_0 + \cdots + u_6 Y_6 = u_1 Y_1 - u_3 Y_3 = u_2 Y_2 - u_3 Y_3 =
u_4 Y_4 - u_6 Y_6 = u_5 Y_5 - u_6 Y_6 = 0;  \]
\[ Y_0^3 = Y_1 Y_2 Y_3 =  Y_4 Y_5 Y_6 \]
has a non-zero solution.

If we put
\[ A = u_3 Y_3,\; B = u_6 Y_6,\; C = u_0 Y_0 \]
then the last system can be rewritten as
\[ 3A + 3B + C =  A^3 + z_1 C^3 = B^3 +
z_2 C^3 = 0.  \]

So  the discriminant the  two-parameter family
is the resultant of two binary homogeneous cubic
equations in $A$ and $B$:
\[ 27z_1 (A + B)^3 - A^3 = 0,\; 27z_2 (A + B)^3 - B^3 = 0. \]

Put $27z_1 = x$, $27z_2 =y$.

\begin{prop}
The discriminant of the $2$-parameter family of Calabi-Yau $3$-folds
$\hat{Z}_f$ is
\[ {\rm Disc}\,f = 1 - (x+y) + 3(x^2 - 7xy + y^2) -
(x^3 + 3x^2y + 3xy^2 + y^3). \]
\end{prop}

\subsection{The diagonal one-parameter subfamily}

\noindent

We consider the diagonal  one-parameter
subfamily of K\"ahler structures $\eta$ on $V$ which are invariant under
the natural involution on $H^{1,1}(V)$, i.e., we assume that
\[ \int_{\gamma_1} \eta = \int_{\gamma_2} \eta. \]
This is equivalent to the substitution $z = z_1 = z_2$.

\begin{rem}
{\rm In this case we obtain the  one-parameter family of mirrors
\[ f_{\psi}(X) = X_1 +  X_2 +  (X_1 X_2)^{-1} +
X_3 +  X_4 +  (X_3 X_4)^{-1} - {3}{\psi}  = 0, \; {\psi}^3 = (27 z)^{-1} \]
which is an analog to mirrors of quintic $3$-folds \cite{cand2}. }
\end{rem}

It is easy to check that the discriminant of $f_{\psi}(X)$ vanishes
exactly when $\psi = \alpha + \beta$, where $\alpha^3 = \beta^3 =1$, i.e.,
$\psi^3 \in \{ 8, -1 \}$, or $z \in \{ -(3)^{-3}, (2\cdot3)^{-3} \}$.

The monodromy invariant period function is
\[ F_0(z) = \Phi_0(z,z) = \sum_{n \geq 0}
\left( \sum_{k + m = n}
\frac{(3n)!}{(k!)^3 (m!)^3} \right) z^n. \]
It satisfies an ordinary Picard-Fuchs differential equation
\[ {\cal D}\;: \; (\Theta^4  +
\sum_{i =0}^3 C_i(z)\Theta^i )F(z)  = 0,\;\; \Theta =
z\frac{\partial}{\partial z}. \]

We compute  the Picard-Fuch differential equation ${\cal E}$ for
$F_0(z)$ is from the  recurrent formula for the coefficients
\[ a_n = \sum_{k + m = n}
\frac{(3n)!}{(k!)^3 (m!)^3} = \frac{(3n)!}{(n!)^3} (
\sum_{k =0}^n { n \choose k }^3) \]
in the power expansion
\[ F_0(z) = \sum_{n \geq 0} a_n z^n .\]

\begin{prop} {\rm (\cite{stienstra}) }
Let
\[ b_n = \sum_{k =0}^n { n \choose k }^3. \]
Then the numbers $b_n$ satisfy the recurrent relation
\[ (n+1)^2 b_{n+1} = (7 n^2 + 7n + 2) b_{n} + 8 n^2 b_{n-1}. \]
\end{prop}

\begin{coro}
The numbers $a_n$ satisfy the recurrent relation
\[ (n+1)^4 a_{n+1} =
3( 7n^2 + 7n + 2) (3n +2) (3n +1) a_n +
72( 3n + 2)( 3n + 1)( 3n - 1)( 3n - 2)a_{n-1}. \]
\end{coro}

\begin{coro}
The monodromy invariant period function $F_0(y)$ is annihilated  by
  the differential operator ${\cal P}$ :
\[ \Theta^4 - 3z (7 \Theta^2 + 7\Theta + 2) ( 3\Theta + 1)
( 3 \Theta + 2) - 72 z^2 ( 3\Theta + 5)
(3 \Theta + 4) ( 3\Theta + 2)( 3\Theta + 1). \]
\end{coro}

The last operator can be rewritten also as
\[  (1 - 216z)(1 + 27z) \Theta^4  -
 54z( 7 + 432z) \Theta^3  - 3 z (10584 z + 95) \Theta^2  -
 48z (351z + 2)\Theta - 12z - 2880z^2 .\]

In particular, one has the coefficient
\[ C_3(z) = \frac{-54z( 7 + 432z)}{(1 - 216z)(1 + 27z)}. \]

The $z$-normalized Yukawa coupling $K_z^{(3)}$ is the solution to the
differential equation
\[ \frac{dK_z^{(3)}}{dz} =  \frac{27( 7 + 432z)}
{(1 - 216z)(1 + 27z)} K_z^{(3)}. \]

Let $H$ be the cohomology class in $H^2(V, {\bf Z})$ such that
$\langle H, \gamma_1 \rangle = \langle H, \gamma_2 \rangle = 1$.
Since $H^3 = 18$, we obtain the normalization condition
\[ K_{z}^{(3)}(0) = 18. \]

Applying the general algorithm in \ref{log-q}, we find the  $q$-expansion of
the
$z$-coordinate
\[ z(q) = q - 48 q^2 - 18q^3 + 7976 q^4 - 1697115 q^5 + O(q^6), \]
and the $q$-expansion of the $q$-normalized Yukawa coupling is
\[ K_q^{(3)} = 18 + 378q + 69498 q^2 + 7724862 q^3 +
1030043898 q^4 + 132082090128 q^5 + O(q^6). \]

We expect that
\[ K_q^{(3)}  = 18 + \sum_{d =1}^{\infty} \frac{n_d d^3 q^d}{1- q^d}. \]
where  $n_d$  are  predictions  for numbers
rational curves of  degree $d$ relative to the ample
divisor of type $(1,1)$ on $V$. In paricular, one has $n_1 = 378.$

\subsection{Lines on a generic Calabi-Yau 3-fold in
${\bf P}^2 \times {\bf P}^2$}

\noindent

We show how to check the prediction for the number  of lines
on a generic Calabi-Yau 3-fold in ${\bf P}^2 \times {\bf P}^2$.

First we formulate one lemma which will be useful in the sequel.

\begin{lem}
Let $M$ be a complete algebraic variety, ${\cal L}_1$ and
${\cal L}_2$ two invertible sheaves on $M$ such that the  projectivizations
${\bf P}({\cal L}_i) = {\bf P} (H^0(M, {\cal L}_i))$ $(i =1,2)$ are nonempty.
Define the morphism
\[ p_{\lambda} \; : \; {\bf P}({\cal L}_1) \times {\bf P}({\cal L}_2)
\rightarrow {\bf P}({\cal L}_1 \otimes {\cal L}_2) = {\bf P} ( H^0(M,
{\cal L}_1 \otimes {\cal L}_2))\]
by the natural mapping
\[ \lambda\; : \; H^0(M, {\cal L}_1) \otimes H^0(M, {\cal L}_2)
\rightarrow H^0(X, {\cal L}_1 \otimes {\cal L}_2). \]
Then the pullback $p_{\lambda}^* {\cal O}(1)$ of the ample generator
${\cal O}(1)$ of
the Picard group of ${\bf P}({\cal L}_1 \otimes {\cal L}_2)$ is isomorphic
to ${\cal O}(1,1)$ on ${\bf P}({\cal L}_1) \times {\bf P}({\cal L}_2)$.
\label{mapping}
\end{lem}

{\bf Proof.} The statement follows immediately from the fact that
$\lambda$ is bilinear.
 \hfill $\Box$

\begin{prop}
A generic Calabi-Yau hypersurface in ${\bf P}^2 \times {\bf P}^2$
contains $378$ lines relative to the ${\cal O}(1,1)$-polarization.
\end{prop}

{\bf Proof. } There are two possibilities for the type  of lines:
$(1,0)$ and $(0,1)$. By symmetry, it is sufficient to consider only
$(1,0)$-lines whose projections on the second factor in ${\bf P}^2 \times
{\bf P}^2$ are points. Let
\[ \pi_2 \;: \; V \rightarrow {\bf P}^2 \]
be the projection of $V$ on the second factor. Then for every point
$p \in {\bf P}^2$ the fiber $\pi_2^{-1}(p)$ is a cubic
in ${\bf P}^2 \times {p}$. We want to calculate the number of
those fibers  $\pi_2^{-1}(p)$ which are unions of a line $L$ and a conic $Q$
in ${\bf P}^2 \times {p}$.
The space of the reducible cubics $L \cup Q$ is isomorphic to the
image $\Lambda \subset {\bf P}^9 = {\bf P}({\cal O}_{{\bf P}^2}(3))$
of the morphism
\[ {\bf P}( {\cal O}_{{\bf P}^2}(1)) \times
{\bf P}( {\cal O}_{{\bf P}^2}(2))) = {\bf P}^2 \times {\bf P}^5
\rightarrow {\bf P}^9 = {\bf P}( {\cal O}_{{\bf P}^2}(3)).  \]
By \ref{mapping}, $\Lambda$ has codimension 2 and degree $21$.

On the other hand, a generic Calabi-Yau hypersurface $V$ defines a
generic Veronese embedding
\[ \phi \; : \; {\bf P}^2
\hookrightarrow {\bf P}^9 = {\bf P}( {\cal O}_{{\bf P}^2}(3)),
\; \; \phi(p) = \pi_2^{-1} (p). \]
The degree of the image $\phi({\bf P}^2)$ is $9$. The number of $(1,0)$-lines
is the intersection number of two subvarieties $\phi({\bf P}^2)$ and
$\Lambda$ in ${\bf P}^9$, i.e., $9 \times 21 = 189$. Thus, the total amount
of lines is $2 \times 189 = 378$.  \hfill $\Box$
\bigskip

\section{Further  examples}

\noindent

In this section we consider more  examples of Calabi-Yau $3$-folds $V$
obtained as complete intersections in product of projective spaces. In all
these examples for simplicity we restrict ourselves to one-parameter
subfamilies invariant under permutations of factors. The latter allows to
apply the Picard-Fuchs operators of order $4$ to the calculation of
predictions for numbers of rational curves on Calabi-Yau $3$-folds with
$h^{1,1} > 1$.
\bigskip

\subsection{Calabi-Yau $3$-folds in
${\bf P}^1 \times {\bf P}^1 \times {\bf P}^1
\times {\bf P}^1$}

\noindent

We consider the diagonal subfamily of K\"ahler classes on Calabi-Yau
hypersurfaces of degree $(1,1,1,1)$ in $({\bf P}^1)^4$. Repeating the same
procedure as for hypersurfaces of degree $(3,3)$ in
${\bf P}^2 \times {\bf P}^2$, we obtain:

\begin{center}

\begin{tabular}{|c|c|} \hline
& \\
$F_0(z)$ & $ {\displaystyle \sum_{n =0}^{\infty} \left(
\sum_{ k_1 + k_2 + k_3 + k_4 = n  }
 \frac{( 2k_1 + 2k_2 + 2k_3 + 2k_4 )! }
{ (k_1!)^2 (k_2!)^2 (k_3!)^2 (k_4!)^2 } \right) z^n   }$ \\
& \\
\hline
& \\
${\cal P}$ & $\Theta^4 - 4 z (5\Theta^2 + 5 \Theta + 2) (2\Theta +1) +
64 z^{2} (2 \Theta +3) (2 \Theta +1) (2\Theta + 2 )^2 $ \\
& \\
\hline
& \\
$ K_z^{(3)}$ & ${\displaystyle  \frac{48}{(1 - 64z)(1 -16z)}  }$ \\
& \\
\hline
& \\
$ K_q^{(3)}$ &  $   48 + 192q + 7872 q^2 + 278400 q^3 + 9445056 q^4 +
315072192 q^5 + O(q^6)$ \\
& \\
\hline
& \\
$ n_i$  & $ n_1 = 192,\; n_2 = 960, \; n_3 = 10304,\; n_4 = 147456,\;
n_5 = 2520576$ \\
& \\
\hline
\end{tabular}

\end{center}

\begin{prop}
The number of  lines on a generic Calabi-Yau hypersurface
in  ${\bf P}^1 \times {\bf P}^1 \times {\bf P}^1
\times {\bf P}^1$ relative to the $(1,1,1,1)$-polarization is
equal to $192$.
\end{prop}

{\bf Proof. } Let $f$ be the polynomial of degree $(2,2,2,2)$ defining
a Calabi-Yau hypersurface $V$ in $({\bf P}^1)^4$. If $V$ contains
a $(0,0,0,1)$-curve whose projection on the product of first three
${\bf P}^1$ is a point $(p_1,p_2,p_3)$, then all three coefficients of
the binary quadric obtained from $f$ by substitution of $(p_1,p_2,p_3)$
must vanish. Hence, the number of $(0,0,0,1)$ curves on $V$ equals the
intersection number of $3$ hypersurfaces of degree $(2,2,2)$ in
${\bf P}^1 \times {\bf P}^1 \times {\bf P}^1$. This number is $48$. By
symmetry, the total amount of lines on $V$ is $4 \times 48 = 192$.  \hfill
$\Box$

\begin{prop}
The number of  conics on a generic Calabi-Yau hypersurface
in  ${\bf P}^1 \times {\bf P}^1 \times {\bf P}^1
\times {\bf P}^1$ with respect to the $(1,1,1,1)$-polarization is
equal to $960$.
\end{prop}

{\bf Proof.}
By symmetry, it is sufficient to compute the number of rational
curves of type $(0,0,1,1)$. Let $M$ be the product of first two
${\bf P}^1$  in $({\bf P}^1)^4$. Then we obtain the natural embedding
\[ \phi \; : \; M \hookrightarrow {\bf P}^8 =
{\bf P}({\cal O}_{{\bf P}^1 \times {\bf P}^1}(2,2)). \]
On the other hand, the points on $M$ corresponding to
projections of $(0,0,1,1)$-curves on $V$ are intersections
of $\phi(M)$ with the $6$-dimensional subvariety $\Lambda \subset {\bf P}^8$
which is the image of the morphism
\[ \phi' \; : \;
{\bf P}({\cal O}_{{\bf P}^1 \times {\bf P}^1}(1,1))^2 = {\bf P}^3 \times
{\bf P}^3 \rightarrow {\bf P}^8 =
{\bf P}({\cal O}_{{\bf P}^1 \times {\bf P}^1}(2,2)). \]
The image $\phi(M)$ has degree $8$. On the other hand,
$\phi$ has degree two onto its image. Hence,  the subvariety $\Lambda$ has
degree $10$. Hence, we obtain $8 \times 20 = 160$ points on $M$.
There are $6$ possibilities
for the choice of the type of conics. Thus, the total amount of
conics is $6 \times 160 = 960$. \hfill $\Box$
\bigskip

\subsection{Complete intersections of three hypersurfaces
in ${\bf P}^2 \times {\bf P}^2
\times {\bf P}^2$}

\noindent

We consider two examples of $3$-dimensional
complete intersections with trivial canonical
class in $({\bf P}^2)^3$.
\medskip

{\bf  Calabi-Yau  complete intersections of
$3$ hypersurfaces of degree $(1,1,1)$:}

\begin{center}

\begin{tabular}{|c|c|} \hline
& \\
$F_0(z)$ & $ {\displaystyle  \sum_{n =0}^{\infty} \left(
\sum_{  k + m + l = n }  \frac{(( k+ m +l )!)^3 }
{ (k!)^3 (m!)^3 (l!)^3 } \right) z^n  }$ \\
& \\
\hline
& \\
${\cal P}$ & $ 25 \Theta^4 -
15 z (5 + 30 \Theta + 72 \Theta^2 + 84 \Theta^3 +
51 \Theta^4)$  \\
& $ + 6z^2 (15 + 155 \Theta + 541 \Theta^2 + 828 \Theta^3 +
531 \Theta^4)$ \\
&  $- 54 z^3 (1170 + 3795 \Theta + 4399 \Theta^2 + 2160 \Theta^3 +
423 \Theta^4)$ \\
& $ + 243z^4 (402 + 1586\Theta + 2270 \Theta^2 + 1386 \Theta^3 +
279 \Theta^4) - 59049 z^5 (\Theta + 1)^4 $ \\
& \\
\hline
& \\
$ K_z^{(3)}$ & ${\displaystyle  \frac{90 + 162 z}{(27z -1)(27z^2 +1)}  }$ \\
& \\
\hline
& \\
$ K_q^{(3)}$ & $90  + 108 q + 2916 q^2 + 57456 q^3 + 834084 q^4 +
13743108 q^5 + O(q^6)$ \\
& \\
\hline
& \\
$n_i$  & $n_1 = 108,\; n_2 = 351, \; n_3 = 2124,\; n_4 = 12987,\;
n_5 = 109944  $ \\
& \\
\hline
\end{tabular}

\end{center}

\begin{prop}
A generic complete intersection of 3 hyper\-sur\-faces of
deg\-ree $(1,1,1)$ in ${\bf P}^2 \times {\bf P}^2
\times {\bf P}^2$ contains $108$ lines relative to the
${\cal O}(1,1,1)$-polarization.
\end{prop}

{\bf Proof.}  Let $V$ be the complete intersection of three
generic hypersurfaces $V_1$, $V_2$, $V_3$ in
$M_1 \times M_2 \times M_3$, where $M_i \cong {\bf P}^2$.

By symmetry, it is sufficient to consider lines having
the class $(0,0,1)$ whose projections on $M_1 \times M_2$ are
 points.  There is the morphism
\[ \phi \; : \; M_1 \times M_2 \rightarrow {\bf P}^8 = {\bf P}(E) \]
where $E$ is the space of all $3\times3$-matrices. By definition,
$\phi$ maps a point $(p_1, p_2) \in M_1 \times M_2$ to the
matrix of coefficents of there linear forms obtained from
the substitution of $p_1$ and $p_2$ in the equations of $V_1, V_2$, and
$V_3$. The morphism $\phi$ is the Segre embedding and its image has
degree $6$. On the other hand, if a point $(p_1,p_2) \in M_1 \times M_2$
is a projection of a $(0,0,1)$-curve on $V$, then the image
$\phi(p_1,p_2)$ must correspond to a matrix of rank $1$ in $E$.
Thus, the number of $(0,0,1)$-curves equals $6\times 6 = 36$, the
intersection number of two Segre subvarieties in ${\bf P}^8$.
So the number of lines on $V$ is $3 \times 36 = 108$. \hfill $\Box$
\bigskip

{\bf Abelian $3$-folds:} The complete intersection of three hyper\-sur\-fa\-ces
of
deg\-rees $(3,0,0)$, $(0,3,0)$, and $(0,0,3)$ are
are abelian $3$-folds constructed by taking products of $3$ elliptic cubic
curves in ${\bf P}^2$. Although abelian
varieities are not Calabi-Yau manifolds from view point of algebraic
geometers, these manifolds also  present interest for physicists.

\begin{center}

\begin{tabular}{|c|c|} \hline
& \\
$F_0(z)$ & $  {\displaystyle \sum_{
p + q + r = n } \left(
 \frac{((3p)!)^3 ((3q)!)^3 ((3r)!)^3}{(p!)^3 (q!)^3 (r!)^3} \right)z^n }$ \\
& \\
\hline
& \\
${\cal P}$ & $\Theta^4 - 3z (6 + 29\Theta + 56 \Theta^2 +
54 \Theta^3 + 27 \Theta^4) $ \\
 & $ + 81 z^2(27 \Theta^2 + 54 \Theta + 40)(\Theta +1)^2 $ \\
& $ - 2187z^3 (3\Theta + 5)(3\Theta +4)(\Theta +2)(\Theta +1) $ \\
& \\
\hline

& \\
$ K_q^{(3)}$ &162  \\
& \\
\hline
\end{tabular}

\end{center}

Thus, we obtain that all Gromov-Witten invariants for the abelian $3$-folds are
zero which agrees with the fact that there are no rational curves
on abelian varieties.
\bigskip

\subsection{Calabi-Yau $3$-folds  in ${\bf P}^3 \times {\bf P}^3$}

\noindent

{\bf Complete intersections of a  hypersurface of
degree $(2,2)$ and $2$ copies of hypersurfaces of degree $(1,1)$: }

\newpage

\begin{center}

\begin{tabular}{|c|c|} \hline
& \\
$F_0(z)$ & ${\displaystyle \sum_{n =0}^{\infty} \left(
\sum_{ k + m =n } \frac{ (2 k + 2m)! ((k +m )!)^2 }
{(k!)^4 (m!)^4 } \right) z^n  }$ \\
& \\
\hline
& \\
${\cal P}$ & $ \Theta^4 - 4z(3\Theta + \Theta +1) (2\Theta +1)^2 -
 4z^2(4\Theta + 5) (4\Theta + 6)(4\Theta + 2)(4\Theta + 3) $ \\
& \\
\hline
& \\
$ K_z^{(3)}$ & ${\displaystyle  \frac{40}{(1 + 16z)(1 - 64z)}  }$ \\
& \\
\hline
& \\
$ K_q^{(3)}$ & $ 40 + 160 q + 12640q^2 + 393280 q^3 + 17420640q^4 +
662416160 q^5 + O(q^6) $ \\
& \\
\hline
& \\
$ n_i$  & $ n_1 = 160,\, n_2 = 1560, \, n_3 = 14560, \, n_4 = 272000, \,
n_5 = 5299328 $ \\
& \\
\hline
\end{tabular}

\end{center}
\bigskip

\begin{prop}
The number of lines on a generic complete intersections of a  hypersurface of
degree $(2,2)$ and $2$ copies of hypersurfaces of degree $(1,1)$ is equal to
$160$.
\label{3hyper}
\end{prop}

{\bf Proof. } Let $W = {\rm Gr}(2,4) \times {\bf P}^3$ be the $7$-dimensional
variety parametrizing all $(1,0)$-lines on ${\bf P}^3 \times {\bf P}^3$.
Let  ${\cal E}$ be the tautological rank-$2$ locally
free sheaf on ${\rm Gr}(2,4)$. We put
$c_1({\cal E}) = c_1$, $c_2({\cal E}) = c_2$, and $h$ be the first
Chern class of the
ample generator $H$ of ${\rm Pic}({\bf P}^3)$.
Let
$S^2({\cal E})$ be the $2$-nd symmetric power
of  ${\cal E}$. By standard arguments, we obtain:

\begin{lem}
The Chern classes $c_1$, $c_2$ generate the cohomology ring of
${\rm Gr}(2,4)$. The elements $1,\,c_1,\,c_2,\,c_1^2,\,c_1c_2,\, c_1^2c_2$
form a ${\bf Z}$-basis of $H^*({\rm Gr}(2,4), {\bf Z})$,  and
one has the following
\[ c_1^4 = c_2^2 = c_1^2,\; c_1^3 = 2c_1 c_2; \]
 \[ c_1(S^2({\cal E})) = 3c_1,\;
 c_2(S^2({\cal E})) = 2c_1^2 + 4 c_2, \;
c_3(S^2({\cal E})) = 4c_1c_2. \]
Moreover, for any invertible sheaf ${\cal L}$ on ${\bf P}^3$, one has
\[ c_1(S^2({\cal E})\otimes{\cal L}) = 3c_1 + 3c_1({\cal L}),\; \]
\[ c_2(S^2({\cal E})\otimes{\cal L}) = 2c_1^2 + 4 c_2^2 +
2c_1({\cal L})( 3c_1) + 3 c_1^2({\cal L}), \; \]
\[ c_3(S^2({\cal E})\otimes{\cal L}) =
 4c_1c_2 + c_1({\cal L})(2 c_1^2 + 4c_2) + c_1^2({\cal L})(3c_1) +
 c_1^3({\cal L}). \]
\label{chern.cl1}
\end{lem}

Then the number of $(1,0)$-lines equals the following product in the
cohomology ring of $W$:
\[ c_2({\cal E}\otimes{\cal O}(H)) \cdot
c_2({\cal E}\otimes{\cal O}(H)) \cdot
c_3(S^2({\cal E})\otimes{\cal O}(2H)) \]
\[ = (h^2 + c_1 h + c_2)^2 (8 h^3 + 3 c_1 \cdot 4h^2 + 2(c_1^2 + 2c_2)\cdot 2h
+ 4 c_1 c_2) \]
\[ = (8 c_1^2 c_2 h^3 + 4(c_1^2 + 2c_2)^2 h^3 + 24 c_1^2 c_2 h^3 +
8 c_2^2 h^3)  \]
\[ = (8 + 4\times 10 + 24 + 8) c_1^2 c_2 h^3 = 80 c_1^2 c_2 h^3. \]
Thus, the number of $(1,0)$-lines is $80$. By symmetry, the total amount of
lines is $160$.  \hfill $\Box$
\bigskip

{\bf Complete intersections of hypersurfaces of degrees
$(1,1)$, $(1,2)$ and
$(2,1)$:}

\begin{center}

\begin{tabular}{|c|c|} \hline
& \\
$F_0(z)$ & ${\displaystyle  \sum_{n =0}^{\infty} \left(
\sum_{ k + m =n }
\frac{ (2 k + m)! (k + 2m)!(k +m )! }
{(k!)^4 (m!)^4 }  \right)z^n }$ \\
& \\
\hline
& \\
${\cal P}$ & $ 529 \Theta^4 - 23z( 92 + 621 \Theta + 1644 \Theta^2 +
2046 \Theta^3 + 921 \Theta^4) $ \\
& $ - z^2(221168 + 1033528 \Theta + 1772673 \Theta^2 + 1328584 \Theta^3 +
380851 \Theta^4) $ \\
& $ - 2z^3 (-27232 + 208932 \Theta + 1028791 \Theta^2 + 1310172 \Theta^3 +
475861 \Theta^4) $ \\
& $ - 68z^4 (-976 - 1664 \Theta + 5139 \Theta^2 + 14020 \Theta^3 +
8873 \Theta^4) $ \\
& $ + 6936z^5 (3\Theta + 4)(3 \Theta + 2)(\Theta +1)^2 $ \\
& \\
\hline
& \\
$ K_z^{(3)}$ & ${\displaystyle  \frac{46 + 68z}{(54z -1)(z^2 - 11 z -1)}}$ \\
& \\
\hline
& \\
$ K_q^{(3)}$ & $46 + 160 q + 9416 q^2 + 251530 q^3 + 9120968 q^4 +
289172660 q^5 + O(q^6)$ \\
& \\
\hline
& \\
$ n_i$  & $ n_1 = 160,\; n_2 = 1157, \; n_3 = 9310,\; n_4 = 142368,\;
n_5 = 2313380 $ \\
& \\
\hline
\end{tabular}

\end{center}
\bigskip

\begin{prop}
The number of lines on a generic complete intersections of  hypersurfaces of
degrees $(2,1)$, $(1,2)$, and $(1,1)$ is equal to
$160$.
\end{prop}

{\bf Proof. } We use the same notations as in \ref{3hyper}. The number of
$(1,0)$-lines equals the the following product in the
cohomology ring of $W$:
\[ c_2({\cal E}\otimes{\cal O}(H)) \cdot
c_2({\cal E}\otimes{\cal O}(2H)) \cdot
c_3(S^2({\cal E})\otimes{\cal O}(H))   \]
\[ = (h^2 + c_1 h + c_2) \cdot (4 h^2 + 2 c_1 h + c_2) \cdot
(h^3 + 3 c_1h^2 + 2( c_1^2 + 2 c_2) h + 4 c_1 c_2) \]
\[ = (24 c_1^2 c_2 + 2 (5 c_2 + 2c_1^2)(2 c_2 + c_1^2) + 9 c_1^2 c_2 +
c_2^2)h^3 \]
\[ = (24 + 2 ( 5 + 10 + 4 + 4 ) + 9 + 1)c_1^2 c_2 h^3 = 80 c_1^2 c_2 h^3. \]
Thus, the number of $(1,0)$-lines equals $80$. By symmetry, the total amount
of lines is $160$. \hfill $\Box$
\bigskip

{\bf Hypersurfaces in product of two
Del Pezzo surfaces of degree $3$:}

A Calabi-Yau hypersurface in product of two Del Pezzo surfaces of degree
$3$ is a complete intersections of $(1,1)$, $(3,0)$ and
$(0,3)$-hypersurfaces in ${\bf P}^3 \times {\bf P}^3$.

\begin{center}

\begin{tabular}{|c|c|} \hline
& \\
$F_0(z)$ & ${\displaystyle \sum_{n =0}^{\infty} \left(
\sum_{  k + m =n } \frac{ (3 k)! (3m)!(k +m )! }
{(k!)^4 (m!)^4 }  \right) z^n  }$ \\
& \\
\hline
& \\
${\cal P}$ & $ \Theta^4 - 3z(4 + 23 \Theta + 53 \Theta^2 +
60 \Theta^3 + 48 \Theta^4)$ \\
& $ + 9z^2(304 + 1344\Theta + 2319\Theta^2 + 1980\Theta^3 + 873 \Theta^4)$ \\
& $ - 162z^3(800 + 3348 \Theta + 5259\Theta^2 + 3888\Theta^3 +
1269 \Theta^4)$ \\
& $ + 2916z^4(688 + 2952\Theta + 4653 \Theta^2 + 3240\Theta^3 + 891 \Theta^4
$\\
& $ - 1417176z^5 (3\Theta +4)(3\Theta + 2)(\Theta +1)^2 $ \\
& \\
\hline
& \\
$ K_z^{(3)}$ & ${\displaystyle  \frac{54 -972z}{(1 - 54z)(1 -27z)^2}  }$ \\
& \\
\hline
& \\
$ K_q^{(3)}$ & $54 + 162 q + 7290 q^2 + 119232 q^3 + 3045114 q^4 +
79845912 q^5 + O(q^6)$ \\
& \\
\hline
& \\
$ n_i$  & $ n_1 = 162,\; n_2 = 891, \; n_3 = 4410,\; n_4 = 47466,\;
n_5 = 638766 $ \\
& \\
\hline
\end{tabular}

\end{center}

\begin{prop}
Let $S_1$, $S_2$ be two Del Pezzo surfaces of degree $3$. Then the number
of lines on a generic Calabi-Yau hypersurface $V$ in $S_1 \times S_2$ is
$162$.
\end{prop}

{\bf Proof.}  If $C$ is a line of type $(1,0)$ on $S_1 \times S_2$, then
$\pi_1(C)$ is one of $27$ lines on $S_1$, and $\pi_2(C)$ is a point on $S_2$.
Let ${\cal O}_{S_i}(-K)$ denotes the anticanonical bundle over $S_i$. Then
the zero set of a generic global section $s$  of
$\pi_1^*{\cal O}_{S_1}(-K) \otimes
\pi_2^*{\cal O}_{S_2}(-K)$ defines  a morphism
\[  \phi\; : \; S_2 \rightarrow {\bf P}^3 = {\bf P}({\cal O}_{S_1}(-K)). \]
On the other hand, for any line $L \in S_1$, one has the linear embedding
\[ \phi'\; : \; {\bf P}({\cal O}_{S_1}(-K-L)) \cong {\bf P}^1 \hookrightarrow
{\bf P}^3 = {\bf P}({\cal O}_{S_1}(-K)). \]
The intersection number of ${\rm Im}\,\phi$ and ${\rm Im}\,\phi'$ in
${\bf P}^3$ equals $3$, i.e., one has exactly $3$ lines $C$
on a generic $V$ such that $\pi_1(C) = L$ and $\pi_2(C)$ is a point
on $S_2$. Thus, there are $3 \times 27 = 81$ lines of type $(1,0)$ on $V$.
By symmetry, the total amount of lines is $162$.  \hfill $\Box$

\begin{prop}
Let $S_1$, $S_2$ be two Del Pezzo surfaces of degree $3$. Then the number
of conics on a generic Calabi-Yau hypersurface $V$ in $S_1 \times S_2$ is
$891$.
\end{prop}

{\bf Proof.} If $C$ is a conic of type $(1,1)$ on   $S_1 \times S_2$, then
$L_1 = \pi_1(C)$ is one of $27$ lines on $S_1$, and $L_2 = \pi_2(C)$ is
one of $27$ lines on $S_2$. On the other hand, for any pair of lines
$L_1 \in S_1$, $L_2 \in S_2$, the intersection of
the product $L_1 \times L_2 \subset S_1 \times S_2$ with $V$ is a conic
of  type $(1,1)$. So we obtain $27 \times 27 = 729$ conics of type $(1,1)$
on $V$. On the other hand, the number of $(2,0)$ and $(0,2)$ conics is
obviously equals to the number of $(1,0)$ and $(0,1)$ lines. Thus, the total
number of conics is equal to $729 + 162 = 891$.  \hfill $\Box$
\bigskip

\subsection{Calabi-Yau $3$-folds in
${\bf P}^4 \times {\bf P}^4$}

\noindent

{\bf  Complete intersection of hypersuraces of degrees
$(2,0)$, $(0,2)$, and $3$ copies of hypersurfaces of degree $(1,1)$
:}

\begin{center}

\begin{tabular}{|c|c|} \hline
& \\
$F_0(z)$ & ${\displaystyle  \sum_{n =0}^{\infty} \left(
\sum_{ k + m  = n }
\frac{(( k+ m )!)^3 (2k)!(2m)!}
{ (k!)^5 (m!)^5  } \right) z^n  }$ \\
& \\
\hline
& \\
${\cal P}$ & $25 \Theta^4 - 20z(5 + 30\Theta + 72 \Theta^2 + 84 \Theta^3 +
36 \Theta^4) $ \\
 & $- 16 z^2 (- 35 - 70 \Theta + 71 \Theta^2 + 268 \Theta^3
 + 181 \Theta^4)  $ \\
 & $+ 256 z^3 (\Theta +1)(165 + 375\Theta + 248\Theta^2 + 37\Theta^3) $ \\
 & $+ 1024z^4 (59 + 232 \Theta + 331\Theta^2 + 198\Theta^3 + 39 \Theta^4) +
 32768 z^5 (\Theta +1)^4 $ \\
& \\
\hline
& \\
$ K_z^{(3)}$ & ${\displaystyle
\frac{80 + 128 z}{(1 + 4z)(1 -4z)(1 -32z)}}$ \\
& \\
\hline
& \\
$ K_q^{(3)}$ & $80 + 128q + 3776q^2 + 65792 q^3 + 1299136 q^4 +
23104128 q^5 + O(q^6)$ \\
& \\
\hline
& \\
$ n_i$  & $n_1 = 128,\; n_2 = 456, \; n_3 = 2432,\; n_4 = 20240,\; n_5 =
184832  $ \\
& \\
\hline
\end{tabular}
\end{center}

\begin{prop}
The number of lines on a generic complete intersections
of hypersuraces of degrees
$(2,0)$, $(0,2)$, and $3$ copies of hypersurfaces of degree $(1,1)$
is equal to $128$.
\label{5hyper}
\end{prop}

{\bf Proof. } Let $W = {\rm Gr}(2,5) \times {\bf P}^4$ be the $10$-dimensional
variety parametrizing all $(1,0)$-lines on ${\bf P}^4 \times {\bf P}^4$.
Let  ${\cal E}$ be the tautological rank-$2$ locally
free sheaf on ${\rm Gr}(2,5)$. We put
$c_1({\cal E}) = c_1$, $c_2({\cal E}) = c_2$, and $h$ be the first
Chern class of the
ample generator $H$ of ${\rm Pic}({\bf P}^4)$.
Let  $S^2({\cal E}))$ be the $2$-nd symmetric power
of  ${\cal E}$. Again, by standard arguments, we obtain:

\begin{lem}
The Chern classes $c_1$, $c_2$ generate the cohomology ring of
${\rm Gr}(2,5)$. The elements $1,\,c_1,\,c_2,\,c_1^2,\,c_1c_2,\,c_1^3,\,
 c_1^2c_2, \, c_1^4, \, c_1c_2^2,\, c_2^3$
form a ${\bf Z}$-basis of $H^*({\rm Gr}(2,5), {\bf Z})$,  and
satisfy the following relations
\[ c_1^4 c_2 = 2 c_1^2 c_2^2 = 2 c_2^3,\; c_1^5 = 5 c_1 c_2^2,\;
c_1^6 = 5c_1^2 c_2^2 = 5c_2^3, \; c_1^3 c_2 = 2 c_1 c_2^2. \]
\end{lem}

Then the number of $(1,0)$-lines equals the following product in the
cohomology ring of $W$:

\[ c_1({\cal O}(H)\cdot (c_2(S^2({\cal E}))^3 \cdot c_3(S^2({\cal E})) \]
\[ = (2h) \cdot (h^2 + c_1 h + c_2)^3 \cdot (4 c_1 c_2) =
64 c_1^2 c_2^2 h^4. \]

Thus, the number of $(1,0)$-lines is $64$. By symmetry, the total amount of
lines is $128$.  \hfill $\Box$
\bigskip

{\bf  Complete intersection of $5$ hypersurfaces
of degree $(1,1)$:}
\begin{center}

\begin{tabular}{|c|c|} \hline
& \\
$F_0(z)$ & ${\displaystyle   \sum_{n =0}^{\infty} \left(
\sum_{ k + m  = n }
\frac{(( k+ m )!)^5 }
{ (k!)^5 (m!)^5  } \right) z^n }$ \\
& \\
\hline
& \\
${\cal P}$ & $49 \Theta^4 - 7z(14 + 91 \Theta + 234 \Theta^2 + 286 \Theta^3 +
155 \Theta^4)$ \\
& $ - z^2(15736 + 66094 \Theta + 102261 \Theta^2
+ 680044 \Theta^3  + 16105 \Theta^4) $ \\
&  $+ 8z^3 (476 + 3759 \Theta +  9071 \Theta^2 + 8589 \Theta^3 +
2625 \Theta^4)$ \\
& $ - 16 z^4 (184 + 806 \Theta + 1439 \Theta^2 + 1266 \Theta^3 + 465 \Theta^4)
+ 512 z^5 (\Theta +1)^4$ \\

& \\
\hline
& \\
$ K_z^{(3)}$ & ${\displaystyle \frac{70 - 40z}{(32z -1)(z^2 - 11 z -1)}
 }$ \\
& \\
\hline
& \\
$ K_q^{(3)}$ & $K_q(q) = 70 + 100 q + 5300 q^2 + 79750 q^3 + 1966900 q^4 +
37143850 q^5 + O(q^6)$ \\
& \\
\hline
& \\
$ n_i$  &  $n_1 = 100,\; n_2 = 650, \; n_3 = 2950,\; n_4 = 30650,\;
n_5 = 297150$  \\
& \\
\hline
\end{tabular}

\end{center}

\begin{prop}
A generic complete intersection of generic $5$ hypersurfaces of degree
$(1,1)$ in ${\bf P}^4 \times {\bf P}^4$ contains $100$ lines.
\end{prop}

{\bf Proof.} We give below two different proofs of the statement.

{\bf I}: We keep the notation from the proof of \ref{5hyper}.
Then the number of $(1,0)$-lines equals the following product in the
cohomology ring of $W$:

\[ (c_2(S^2({\cal E}))^5  =  (h^2 + c_1 h + c_2)^5   \]
\[ = (c_1 h + c_2)65 + 5 (c_1 h + c_2)64 h^2 + 10 (c_1 h + c_2)^3 h^4  =
5 c_1^4 c_2 h^4 + 5 { 4 \choose 2} c_1^2 c_2^2 h^4 + 10 c_2^3 h^4 \]
\[ = ( 5 \times 2 + 5 \times 6 + 10) c_1 ^4 c_2 h^4 = 50  c_1^2 c_2^2 h^4. \]

Thus, the number of $(1,0)$-lines is $50$. By symmetry, the total amount of
lines is $100$.

\medskip

{\bf II}: Let $M = M_1 \times M_2$ where $M_i \cong {\bf P}^4$
($i =1,2$). By symmetry, we consider only lines of type
$(0,1)$ whose projections on $M_1$ are points.
The substitution of a point $p \in M_1$
in the equations of the hypersurfaces $H_1, \ldots, H_5 \subset M$ gives
$5$ linear forms $f_1, \ldots, f_5$ in homogeneous coordinates on
$M_2$. A point $p \in M_1$ is a projection of a $(0,1)$-line on
$H_1 \cap \cdots \cap H_5$ if the system of linear forms
has rank $3$. The space of 5 copies of linear forms can be identified
with the space $L$ of matrices $5\times5$. We are interested in the
determinantal subvariety $D$ in ${\bf P}^{24}$ consisting of matrices of
rank $\leq 3$. The subvariety $D$ has the codimension $4$, the ideal of
$D$ is generated by all $4\times 4$ minors. Using the free graded
resolution of the homogeneous coordinate ring of $D$ as a module over
the homogeneous coordinate ring of ${\bf P}^{24}$, we can compute the
degree of $D$ which is equal to $50$
(The Hilbert-Poincare series equals $( 1 + 4t + 10 t^2 + 20 t^3 +
10 t^4 + 4 t^5 + t^6)/(1 - t)^{21}$). On the other hand, the equatons
of generic hypersurfaces $H_1, \ldots, H_5$ define a generic embedding
\[ {\bf P}^4 \hookrightarrow {\bf P}^{24} \]
of ${\bf P}^4$ as a linear subspace. So the number of lines of type
$(0,1)$ on a generic complete intersection is $50$. Thus, the total
number of lines is $100$. \hfill $\Box$
\bigskip

{\bf Hypersurfaces in product of two
Del Pezzo surfaces of degree $4$:}

A Calabi-Yau hypersurface in the product of two Del Pezzo surfaces of degree
$4$ is a complete intersection of $5$ hypersurfaces in
${\bf P}^4 \times {\bf P}^4$: two copies of type $(2,0)$, two copies of
type $(0,2)$, and one copy of type $(1,1)$.

\begin{center}

\begin{tabular}{|c|c|} \hline
& \\
$F_0(z)$ & ${\displaystyle \sum_{n =0}^{\infty} \left(
\sum_{ k + m =n }
\frac{ ((2 k)!)^2 ((2m)!)^2(k +m )! }
{(k!)^5 (m!)^5 }  \right)z^n  }$ \\
& \\
\hline
& \\
${\cal P}$ & $ 9 \Theta^4 - 4z (6 + 33\Theta + 73 \Theta^2 +
80 \Theta^3 + 64 \Theta^4) $\\
& $ + 128z^2 (75 + 315\Theta + 527 \Theta^2 + 440 \Theta^3 + 194 \Theta^4)$ \\
& $ - 4096z^3 (66 + 261 \Theta + 397 \Theta^2 + 288 \Theta^3 + 94 \Theta^4)$ \\
& $ + 131072 z^4(19 + 77 \Theta + 117 \Theta^2 + 80 \Theta^3 + 22 \Theta^4) -
8388608 z^5 (\Theta + 1)^4 $ \\
& \\
\hline
& \\
$ K_z^{(3)}$ & ${\displaystyle  \frac{96 - 1024z}{(1 - 32z)(1 -16z)^2}  }$ \\
& \\
\hline
& \\
$ K_q^{(3)}$ & $96 + 128 q + 3456 q^2 + 38144q^3 + 572800 q^4 +
9344128 q^5 + O(q^6)$ \\
& \\
\hline
& \\
$ n_i$  & $n_1 = 128,\; n_2 = 416, \; n_3 = 1408,\; n_4 = 8896,\;
n_5 = 74752  $ \\
& \\
\hline
\end{tabular}

\end{center}

\begin{prop}
Let $S_1$, $S_2$ be two Del Pezzo surfaces of degree $4$. Then the number
of lines on a generic Calabi-Yau hypersurface $V$ in $S_1 \times S_2$ is
$128$.
\end{prop}

{\bf Proof.}  If $C$ is a line of type $(1,0)$ on $S_1 \times S_2$, then
$\pi_1(C)$ is one of $16$ lines on $S_1$, and $\pi_2(C)$ is a point on $S_2$.
Let ${\cal O}_{S_i}(-K)$ denotes the anticanonical bundle over $S_i$. Then
the zero set of a generic global section $s$  of
$\pi_1^*{\cal O}_{S_1}(-K) \otimes
\pi_2^*{\cal O}_{S_2}(-K)$ defines  a morphism
\[  \phi\; : \; S_2 \rightarrow {\bf P}^4 = {\bf P}({\cal O}_{S_1}(-K)). \]
On the other hand, for any line $L \in S_1$, one has the linear embedding
\[ \phi'\; : \; {\bf P}({\cal O}_{S_1}(-K-L)) \cong {\bf P}^2 \hookrightarrow
{\bf P}^4 = {\bf P}({\cal O}_{S_1}(-K)). \]
The intersection number of ${\rm Im}\,\phi$ and ${\rm Im}\,\phi'$ in
${\bf P}^3$ equals $4$, i.e., one has exactly $4$ lines $C$
on a generic $V$ such that $\pi_1(C) = L$ and $\pi_2(C)$ is a point
on $S_2$. Thus, there are $4 \times 16 = 64$ lines of type $(1,0)$ on $V$.
By symmetry, the total amount  of lines is $128$.  \hfill $\Box$

\begin{prop}
Let $S_1$, $S_2$ be two Del Pezzo surfaces of degree $3$. Then the number
of conics on a generic Calabi-Yau hypersurface $V$ in $S_1 \times S_2$ is
$416$.
\end{prop}

{\bf Proof.} If $C$ is a conic of type $(1,1)$ on   $S_1 \times S_2$, then
$L_1 = \pi_1(C)$ is one of $16$ lines on $S_1$, and $L_2 = \pi_2(C)$ is
one of $16$ lines on $S_2$. On the other hand, for any pair of lines
$L_1 \in S_1$, $L_2 \in S_2$, the intersection of
the product $L_1 \times L_2 \subset S_1 \times S_2$ with $V$ is a conic
of  type $(1,1)$. So we obtain $16 \times 16 = 256$ conics of type $(1,1)$
on $V$.

In order to compute the number of $(2,0)$-conics, we notice that
$S_1$ has  exactly $10$  of conic bundle structures. Moreover, these
conic bundle structures can be divided into $5$ pairs such that
the union of degenerate fibers of each pair is the set of all
$16$ lines on $S_1$.
A  generic global section $s$  of
$\pi_1^*{\cal O}_{S_1}(-K) \otimes
\pi_2^*{\cal O}_{S_2}(-K)$ defines the anticanonical embedding
\[  \phi\; : \; S_2 \hookrightarrow {\bf P}^4 = {\bf P}({\cal O}_{S_1}(-K)). \]
On the other hand, the points $p \in S_2$ such that $\phi(p)$ splits into
the union of two conics $C_1 \cup C_2$ are exactly intersection points
of $\phi(S_2)$ and the image of the embedding
\[ \phi' \; : \; {\bf P}({\cal O}_{S_1}(C_1) \times {\cal O}_{S_1}(C_2)
\cong {\bf P}^1 \times {\bf P}^1 \hookrightarrow {\bf P}^4 =
{\bf P}({\cal O}_{S_1}(-K)). \]
Since the image of $\phi'$ has degree $2$, we obtain $8$ points $p \in S_2$.
Each such a point yields $2$ conics on $\pi^{-1}_2(p)$.
Therefore, for each of $5$ pairs of conic bundle structures we have
$16$ $(2,0)$-conics.

Thus, the total
number of conics is equal to $256 + 2 \times 80 = 416$.
\hfill $\Box$


\bigskip



\begin{thebibliography}{99}

\bibitem{batyrev.dual} V.V. Batyrev, {\em Dual Polyhedra and
Mirror Symmetry for Calabi-Yau Hypersurfaces in Toric Varieties},
Nov. 18 (1992) Preprint, Uni-Essen.

\bibitem{batyrev.var} V. V. Batyrev, {\em Variations of the Mixed Hodge
Structure of Affine Hypersurfaces in Algebraic Tori},
Duke Math. J. {\bf 69}, (1993), 349-409.

\bibitem{batyrev.quant} V. V. Batyrev, {\em Quantum Cohomology
Rings of Toric Manifolds}, Preprint MSRI, (1993).

\bibitem{cand0} P. Candelas, A. M. Dale, C .A. L\"utken, R. Schimmrigk,
{\em Complete Intersection Calabi-Yau Manifolds I, II}, Nucl. Phys.
{\bf B 298} (1988) 493-525, Nucl. Phys. {\bf B 306} (1988) 113-136.

\bibitem{cand01} P. Candelas, M. Lynker, R. Schimmrigk, {\em Calabi-Yau
Manifolds in Weighted ${\bf P}^4$}, Nucl. Phys., {\bf B 341} (1990), 383-402.


\bibitem{cand02} P. Candelas, A. He, {\em On the Number of Complete
Intersection
of Calabi-Yau Manifolds},  Commun. Math. Phys. {\bf 135} (1990) 193-199.

\bibitem{cand03} P. Candelas and X. de la Ossa, Nucl. Phys. {\bf 342} (1990)
246.

\bibitem{cand2} P. Candelas, X.C. de la Ossa, P.S. Green, and L. Parkes,
{\em A pair of Calabi-Yau manifolds as an exactly soluble
superconformal theory}, Nuclear Phys. B {\bf 359} (1991), 21-74.

\bibitem{cere} A. Ceresole, R. D'Auria, S. Ferrara, W. Lerche and
J. Louis, {\em Picard-Fuchs Equations and Special Geometry},
CERN-TH.6441/92, UCLA/92/TEP/8, CALT-8-1776, POLFIS-TH.8/92.

\bibitem{deligne1} P. Deligne, {\em Letter to D. Morrison},
6 November 1991.

\bibitem{dixon} L. Dixon, in {\em Superstrings. Unified Theories and
Cosmology 1987}, G. Furlan et al., eds., World Scientific 1988.

\bibitem{ferrara.louis} S. Ferrara, J. Louis.  Phys. Lett. {\bf B 273},
(1990) 246.

\bibitem{font} A. Font, {\em Periods and duality symmetries in
Calabi-Yau com\- pac\-ti\-fi\-ca\-tions}, Nucl. Phys.
{\bf B 389} (1993) 153


\bibitem{gel1} I.M. Gelfand, M.M. Kapranov, and A.V. Zelevinsky,
{\em Hypergeometric functions and toric manifolds},  Functional  Anal.
Appl. {\bf 28:2} (1989), 94-106.

\bibitem{gelfand} I.M. Gelfand, M.M. Kapranov, and A.V. Zelevinsky
{\em Discriminants of polynomials in several variables and triangulations
of Newton polytopes}, {\sl Algebra i analiz (Leningrad Math. J.)}
{\bf 2} (1990) 1-62.


\bibitem{greene} B.R. Greene and M.R. Plesser, {\em Duality in Calabi-Yau
moduli space}, Nuclear Phys. B {\bf 338} (1990), 15-37.

\bibitem{greene1} B.R. Greene, D.R. Morrison, M.R. Plesser, {\em Mirror
manifolds in higher dimension}, in preparation.

\bibitem{katz1} S. Katz, {\em Rational curves on Calabi-Yau $3$-folds}, in
{\sl Essays on Mirror Manifolds} (Ed. S.-T. Yau), Int. Press Co., Hong Kong,
(1992) 168-180.


\bibitem{katz2} S. Katz, {\em Rational curves on Calabi-Yau manifolds:
verifying of predictions of Mirror Symmetry},  Preprint,
alg-geom/9301006 (1993)


\bibitem{klemm1} A. Klemm, S. Theisen, {\em Consideration of One Modulus
Calabi-Yau Compactification$:$ Picard-Fuchs Equation, K\"ahler Potentials
and Mirror Maps},  Nucl. Phys. {\bf B 389} (1993) 153.

\bibitem{klemm2} A. Klemm, S. Theisen, {\em Mirror Maps and Instanton Sums for
Intersections in Weighed Projective Space}, LMU-TPW 93-08, Preprint
 April 1993.

\bibitem{kim} J. K. Kim, C.J. Park and Y. Yoon, {\em Calabi-Yau manifolds
from complete intersections in products of weighted projective spaces},
Physics Letters B, {\bf 224}, (1989), 108-114.

\bibitem{kreu} M. Kreuzer, R. Schimmrigk, H. Skarke, {\em Abelian
Landau-Ginzburg Orbifolds and Mirror Symmetry}, Nucl. Hys. {\bf B 372}
(1992) 61-86.

\bibitem{lerche} W. Lerche, C. Vafa and N. Warner, Nucl. Phys. {\bf B324}
(1989) 427.

\bibitem{lib.teit} A. Libgober, and J. Teitelbaum, {\em Lines on Calabi-Yau
complete intersections, mirror symmetry, and Picard-Fuchs equations},
Duke Math. J., Int. Math. Res Notices {\bf 1} (1993) 29.

\bibitem{lynker} M. Lynker, R. Schimmrighk, Phys. Lett. {\bf 249} (1990) 237.

\bibitem{hyp.geom} A. M. Mathai, R.K. Saxens, {\em Generalized Hypergeometric
Functions with Applications in Statistic and Physical Sciences}, Lect. Notes
Math, {\bf 348}, (1973),

\bibitem{morrison.mirror} D. Morrison, {\em Mirror symmetry and rational
curves on quintic threefolds: a guide for mathematicians},  J. Amer. Math.
Soc. {\bf 6} (1993), 223-247.

\bibitem{morrison.picard} D. Morrison, {\em Picard-Fuchs equations and mirror
maps
for hypersurfaces}, in: {\sl Essay on Mirror Manifolds} (Ed. S.-T.
Yau), Int. Press. Co., Hong Kong, (1992) 241-264.

\bibitem{morrison.comp} D. Morrison, {\em Compactifications of moduli
spaces inspired by mirror symmetry}, Preprint (1993)

\bibitem{morrison.hodge} D. Morrison, {\em Hodge-theoretic aspects of mirror
symmetry}, in preparation


\bibitem{schimm} R. Schimmrigk, {\em A new construction of a three-generation
Calabi-Yau manifold}, Physics Letters B, {\bf 193} (1987), 175-180.

\bibitem{slater} L.J. Slater, {\em Generalized Hypergeometric Functions},
Cambridge University Press, XIII, 1966.

\bibitem{stienstra} J. Stienstra and F. Beukers, {\em On the Picard-Fuchs
Equation and the Formal Brauer Group of Certain Elliptic $K3$-surfaces},
Math. Ann. {\bf 271} (1985), 269-304.

\bibitem{witten0} E. Witten, {\em Topological sigma models}, Commun. Math.
Phys., {\bf 118} (1988) 411-449.

\bibitem{witten} E. Witten, {\em Two-dimensional gravity and intersection
theory on moduli space}, Surveys in Diff. Geom. {\bf 1} (1991) 243-310.

\bibitem{witten2} E. Witten, {\em Mirror Manifolds and Topological
Field Theory}, in
{\sl Essays on Mirror Manifolds} (Ed. S.-T. Yau), Int. Press Co., Hong Kong,
(1992) 120-180.

\end{thebibliography}
\end{document}